\documentclass[11pt]{article}

\usepackage{graphicx,mathtools}
\usepackage{tikz,enumerate, multicol}
\usepackage{float,bbm}
\usepackage{scalerel}
\usepackage{amssymb}
\usepackage{amsfonts}
\usepackage{latexsym}
\usepackage{epsfig}
\usepackage{amsmath}
\usepackage{amsthm}
\usepackage{color}
\usepackage{graphicx}
\usepackage{subfig}
\usepackage{grffile}
\usepackage{tikz-network}
\usetikzlibrary{graphs}

\usetikzlibrary{automata,positioning, math}
\usetikzlibrary{shapes.geometric,arrows,shapes.misc, arrows.meta}
\tikzset{
cnd/.style={
draw, circle, minimum size=10pt, inner sep=0pt, font=\tiny
}}

\newtheorem{ittheorem}{Theorem}
\newtheorem{itlemma}{Lemma}
\newtheorem{itproposition}{Proposition}
\newtheorem{itdefinition}{Definition}
\newtheorem{itremark}{Remark}
\newtheorem{itconjecture}{Conjecture}

\newenvironment{theorem}{\addtocounter{counter}{1}
\begin{ittheorem}}{\end{ittheorem}}

\newenvironment{lemma}{\addtocounter{counter}{1}
\begin{itlemma}}{\end{itlemma}}

\newenvironment{proposition}{\addtocounter{counter}{1}
\begin{itproposition}}{\end{itproposition}}

\newenvironment{definition}{\addtocounter{counter}{1}
\begin{itdefinition}}{\end{itdefinition}}

\newenvironment{remark}{\addtocounter{counter}{1}
\begin{itremark}}{\end{itremark}}

\newenvironment{conjecture}{\addtocounter{counter}{1}
\begin{itconjecture}}{\end{itconjecture}}

\newcommand{\be}[1]{\begin{equation}\label{#1}}
\newcommand{\ee}{\end{equation}}

\newcommand{\bl}[1]{\begin{lemma}\label{#1}}
\newcommand{\el}{\end{lemma}}

\newcommand{\br}[1]{\begin{remark}\label{#1}}
\newcommand{\er}{\end{remark}}

\newcommand{\bt}[1]{\begin{theorem}\label{#1}}
\newcommand{\et}{\end{theorem}}

\newcommand{\bd}[1]{\begin{definition}\label{#1}}
\newcommand{\ed}{\end{definition}}

\newcommand{\bc}[1]{\begin{conjecture}\label{#1}}
  \newcommand{\ec}{\end{conjecture}}

\makeatletter
 \@addtoreset{equation}{section}
 \makeatother

 \makeatletter
 \@addtoreset{counter}{section}
 \makeatother

 \newcounter{counter}[section]

\newcommand{\R}{\mathbb{R}}
\newcommand{\N}{\mathbb{N}}

\newcommand{\ddd}{\mathrm{d}}
\newcommand{\eee}{\mathrm{e}}

\newcommand{\PP}{\mathbb{P}}
\newcommand{\EE}{\mathbb{E}}

\newcommand{\W}{\mathcal{W}}

\newcommand{\AT}{\mathrm{AT}}
\newcommand{\IT}{\mathrm{IT}}
\newcommand{\TP}{\mathrm{TP}}
\newcommand{\C}{\mathrm{C}}
\newcommand{\T}{\mathrm{T}}
\newcommand{\NB}{\mathrm{NB}}
\newcommand{\OB}{\mathrm{OB}}

\newcommand{\cm}{\mathop{\text{CM}_n(\mathbf{d})}}

\definecolor{vrdg}{RGB}{0,0,0}
\newcommand{\vrdg}[1]{{\color{vrdg}{#1}}}
\definecolor{vry}{RGB}{253, 231, 37}

\definecolor{vrg}{RGB}{94,201,98}

\definecolor{vrb}{RGB}{59,82,139}

\definecolor{vrp}{RGB}{68,1,84}

\definecolor{vro}{RGB}{249,142,9}

\definecolor{vrr}{RGB}{188,55,84}

\definecolor{vrnb}{RGB}{13,8,135}

\usepackage[colorlinks=true]{hyperref}
\hypersetup{pdftex, colorlinks=true, linkcolor=vrr, citecolor=vrr,  filecolor=blue,urlcolor=vrb}

\begin{document}


\title{The Triangle Friendship Paradox}

\author{\renewcommand{\thefootnote}{\arabic{footnote}}
Bishakh Bhattacharya
\footnotemark[1]
\\
\renewcommand{\thefootnote}{\arabic{footnote}}
Nitya Gadhiwala
\footnotemark[2]
\\
\renewcommand{\thefootnote}{\arabic{footnote}}
Frank den Hollander
\footnotemark[3]
\\
\renewcommand{\thefootnote}{\arabic{footnote}}
Pradeeptha R Jain
\footnotemark[4]
\\
\renewcommand{\thefootnote}{\arabic{footnote}}
Tejashree Subramanya
\footnotemark[5]
}


\footnotetext[1]{
Indian Statistical Institute, Kolkata, India.\\
email: bishakh.rik@gmail.com
}

\footnotetext[2]{
Department of Mathematics, University of British Columbia, Vancouver, Canada.\\
email: ngadhiwala@math.ubc.ca
}

\footnotetext[3]{
Mathematical Institute, Leiden University, The Netherlands.\\
email: denholla@math.leidenuniv.nl
}

\footnotetext[4]{
International Centre for Theoretical Sciences, Bengaluru, India.\\
email: pradeeptha.jain@icts.res.in
}

\footnotetext[5]{
Indian Institute of Science, Bengaluru, India.\\
email: tejashrees.iisc@gmail.com
}

\maketitle

\begin{abstract}
We consider the generalised friendship paradox, focussing on the number of triangles at a vertex as the relevant attribute. We show that, contrary to the setting where the attribute is the number of edges at a vertex or the number of wedges at a vertex, the average friendship-bias of the number of triangles at a vertex is not always non-negative. We identify classes of finite deterministic graphs for which the bias is non-negative, and provide examples of finite deterministic graphs for which it is not. For certain classes of sparse and dense random graphs, we compute the scaling of the bias in the limit as the number of vertices tends to infinity.

\medskip\noindent
{\it AMS} 2020 {\it subject classifications.}
05C80, 
60C05, 
60F15. 


\medskip\noindent
{\it Acknowledgment.} \vrdg{The authors thank the organisers of the BIRS--CMI research school and workshop in December 2024 at the Chennai Mathematical Institute (CMI) in India, where the research in this paper started: Louigi Addario-Berry, Siva Athreya, Shankar Bhamidi, Serte Donderwinkel and Soumik Pal. The authors also thank BIRS, CMI, ICTS and NBHM for support and hospitality. NG, FdH and PJ thank Siva Athreya for discussions and feedback on the manuscript, and BB thanks Rajat Hazra. FdH is grateful to Rajat Hazra, Vera Meerpoel and Azadeh Parvaneh for discussions on the Generalised Friendship Paradox and to Remco van der Hofstad for discussions on triangles in random graphs.

FdH was supported by the Netherlands Organisation for Scientific Research (NWO) through Gravitation-grant NETWORKS-024.002.003, and by the National Science Foundation (NSF) under Grant No.\ DMS-1928930 while in residence at the Simons Laufer Mathematical Sciences Institute in Berkeley, California, USA during the Spring 2025 semester. PJ acknowledges the grant at ICTS from the Department of Atomic Energy, Government of India under Project No. RTI4001. TS acknowledges the CISCO Fellowship for CNI supported by the Centre for Networked Intelligence, IISc.}
\end{abstract}


\section{Introduction and main results}
\label{s.intro}

The Friendship Paradox (FP) in social networks says that ``on average our friends have more friends than we do'' \cite{F1991}. This observation may at first sight seem paradoxical, yet it is a true fact. \vrdg{The FP has found applications in epidemiology \cite{CF2010}, polling predictions for elections \cite{NK2021}, and various other domains. In \cite{CKN2021} and \cite{PYNSBN2019}, the FP is studied in the context of different random graph models. The FP is considered from a probabilistic point view in \cite{CR2016}, while \cite{HdHP2023} explores the FP in sparse random graphs by using the framework of local weak convergence. 

The Generalized Friendship Paradox (GFP), as the name suggests, is a generalisation of friendship paradox to attributes of the graph other than the degree. Also the GFP has been studied widely in applications. For instance, see \cite{EJ2014}, where the GFP for attributes like collaborations, publications, and citations in co-authorship networks is studied, and \cite{HKL2021}, where the number of followers on social networking sites is studied. We refer the reader to \cite{HdHP2023} for further motivation and a detailed review of the literature for both the FP and the GFP.}

The goal of the present paper is to analyse the GFP for the special case where the attribute assigned to a vertex is the number of triangles the vertex is part of. For certain classes of graphs, both deterministic and random, we compute the sign and the size of the average triangle friendship-bias, thereby providing a partial classification and quantification. We also offer insight into why in this context triangles are harder to deal with than edges or wedges.

\vrdg{The outline of the remainder of this introduction is as follows. In Section~\ref{ss.backmot} we introduce the FP, in Section~\ref{ss.model} the GFP. Our main results are stated in Section~\ref{ss.theorems}. In Section~\ref{sss.star}, we define partially completed star-graphs (see Definition~\ref{def:pstar}) and formulate our two main results, stated in Theorem~\ref{thm:stargraphs} and Theorem~\ref{thm:stargraphsglued}, establishing the triangle friendship-bias for such graphs and its preservation under gluing two of them together. In Section~\ref{sss.sparse}, we discuss two models of sparse random graphs, namely, the Erd\H{o}s-R\'enyi random graph and the Configuration Model, and our main results for these are stated in Theorem~\ref{thm:ER} and Theorem~\ref{thm:CM}. In Section~\ref{sss.dense}, we consider dense random graphs, and Theorem~\ref{thm:densefail} provides conditions under which the triangle friendship paradox holds for the limit or not. We conclude this introduction with Section~\ref{ss.disc}, which discusses the main theorems and formulates a few open problems.}


\subsection{Friendship Paradox}
\label{ss.backmot}

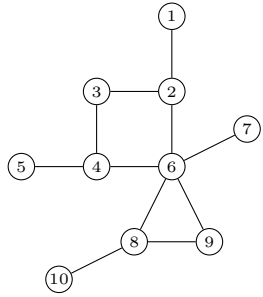
\begin{figure}[htbp]
\centering
\begin{tikzpicture}[scale = 0.25, every node/.style={circle, draw, inner sep = 0pt, minimum size=10pt, font = \tiny}]
    \node (1) at (4,12) {1};
    \node (2) at (4,8) {2};
    \node (3) at (0,8) {3};
    \node (4) at (0,4) {4};
    \node (5) at (-4,4) {5};
    \node (6) at (4,4) {6};
    \node (7) at (8,6) {7};
    \node (8) at (2,0) {8};
    \node (9) at (6,0) {9};
    \node (10) at (-2,-2) {10};
    \draw[-] (1) -- (2);
    \draw[-] (2) -- (3);
    \draw[-] (3) -- (4);
    \draw[-] (4) -- (6);
    \draw[-] (6) -- (2);
    \draw[-] (5) -- (4);
    \draw[-] (7) -- (6);
    \draw[-] (6) -- (8);
    \draw[-] (8) -- (10);
    \draw[-] (6) -- (9);
    \draw[-] (8) -- (9);
\end{tikzpicture}
\caption{\small An example of a finite undirected simple connected graph.}
\label{fig:graph}
\end{figure}

Social networks are modelled as graphs. For $n\in\N$, let $G_n$ be a \emph{finite undirected graph} with $n$ vertices labeled by $[n]=\{1,\ldots,n\}$ (see Figure~\ref{fig:graph}). Let $d_i$ be the degree of vertex $i$. Let $\Delta_{i,n}$ be the \emph{friendship-bias} of vertex $i$, defined as the difference between the average degree of the neighbours of $i$ and the degree of $i$ itself, i.e.,
\[
 \Delta_{i,n} = \bigg[\dfrac{\sum_{j \in [n]} A_{ij} d_j}{d_i}-d_i\bigg] \mathbbm{1}_{\{d_i \neq 0\}}, \qquad i \in [n],
\]
where $(A_{ij})_{i,j \in [n]}$ is the \emph{adjacency matrix} of $G_n$, i.e., $A_{ij}$ is the number of edges between $i$ and $j$ when $i \neq j$ and $A_{ii}$ is twice the number of self-loops at $i$. The FP says that for any choice of $G_n$ the \emph{average friendship-bias} is non-negative, i.e.,
\[
\Delta_{[n]} = \frac{1}{n} \sum_{i \in [n]} \Delta_{i,n} \geq 0.
\]
The proof of this inequality uses the relation $d_i = \sum_{j \in [n]} A_{ij}$ in combination with a symmetrisation argument:
\[
\begin{aligned}
\Delta_{[n]}
&= \dfrac{1}{n} \sum_{ {i \in [n]} \atop {d_i \neq 0} } \bigg(\sum_{j \in [n]}
\frac{A_{ij} d_j}{d_i}-d_i\bigg)
= \dfrac{1}{n} \sum_{ {i \in [n]} \atop {d_i \neq 0} } \sum_{ { j \in [n]} \atop {d_j \neq 0}} A_{ij}
\left(\dfrac{d_j}{d_i}-1\right)\\
&= \dfrac{1}{2n} \sum_{ {i \in [n]} \atop {d_i \neq 0} } \sum_{ { j \in [n]} \atop {d_j \neq 0}} A_{ij}
\left(\dfrac{d_j}{d_i} - 1 + \dfrac{d_i}{d_j} - 1\right)
= \dfrac{1}{2n} \sum_{ {i \in [n]} \atop {d_i \neq 0} } \sum_{ { j \in [n]} \atop {d_j \neq 0}} A_{ij}
\left(\sqrt{\dfrac{d_j}{d_i}}-\sqrt{\dfrac{d_i}{d_j}}\,\,\right)^{2} \geq 0.
\end{aligned}
\]
Equality holds if and only if $i \mapsto d_i$ is constant on each connected component, i.e., all the connected components of $G_n$ are \emph{regular}.

An important question is whether the FP can be \emph{quantified}, and whether $\Delta_{[n]}$ can be analysed for key examples of large random graphs (which are relevant because social networks are typically vast in size and complex in shape). To do so, the key object to look at is the \emph{friendship-bias empirical distribution}
\[
\mu_{n} = \frac{1}{n} \sum_{i \in [n]} \delta_{\Delta_{i,n}}.
\]
Only few papers have addressed the FP from a mathematical viewpoint \cite{PYNSBN2019,CKN2021,HdHP2023}. We refer to the latter reference for a brief overview and relevant literature.


\subsection{Generalised Friendship Paradox}
\label{ss.model}

\vrdg{We first introduce the notion of Generalised Friendship Paradox (GFP), which considers general attributes of vertices that are not necessarily their degrees.}

\begin{definition}
{\rm \vrdg{Let $G_n$ be a finite undirected graph with $n$ vertices labeled by $[n] = \{1,\ldots,n\}$.} Let $x = (x_1,\ldots, x_n)$, where $x_i$ is the \emph{attribute} of vertex $i\in [n]$. In the sequel we focus on \emph{deterministic} attributes, i.e., local functionals of $G$. The \emph{generalised friendship-bias} is defined to be the difference between the average of the attributes of the neighbours of $i$ and the attribute of $i$ itself, i.e.,
\[
\Delta_{i,n}^x = \left[\sum_{j \in [n]} \frac{A_{ij} x_j}{d_i} - x_i\right] 1_{\{d_i \neq 0\}}.
\]
The \emph{Generalised Friendship Paradox} (GFP) is said to hold when
\[
\Delta_{[n]}^x = \frac{1}{n} \sum_{i \in [n]} \Delta_{i,n}^x \geq 0.
\]
}
\hfill$\spadesuit$

\end{definition}

Unlike the FP, the GFP does not always hold. We give a necessary and sufficient criterion. In what follows we will ignore all $i \in [n]$ with $d_i=0$. In the computations below this will be harmless, because $d_i=0$ implies that $A_{ij}=0$ for all $j\in [n]$. Define
\[
\kappa_i = \sum_{j \in [n]} \frac{A_{ij}}{d_j}.
\]
Compute
\[
\frac{1}{n} \sum_{i \in [n]} \kappa_i = \frac{1}{n} \sum_{i \in [n]} \sum_{j \in [n]} \frac{A_{ij}}{d_j}
= \frac{1}{n} \sum_{j \in [n]} \frac{1}{d_j} \left[\sum_{i \in [n]} A_{ij}\right]
=  \frac{1}{n} \sum_{j \in [n]} \frac{1}{d_j} \left[\sum_{i \in [n]} A_{ji}\right]
=  \frac{1}{n} \sum_{j \in [n]} 1 = 1,
\]
i.e., $\frac{1}{n}\kappa$ with $\kappa = (\kappa_i)_{i \in [n]}$ is a probability vector.

\begin{proposition}
$\Delta_{[n]}^x = \mathrm{Cov}(x_U,\kappa_U)$, where $U$ is a vertex that is drawn uniformly at random from $[n]$. Consequently, the GFP holds if and only if $x_U$ and $\kappa_U$ are positively correlated.
\end{proposition}

\begin{proof}
Since
\[
\Delta_{[n]}^x = \frac{1}{n} \sum_{i \in [n]} \left[\,\sum_{j \in [n]} \frac{A_{ij} x_j}{d_i} - x_i\right]
\]
and
\[
\sum_{i \in [n]} \frac{A_{ij}}{d_i} = \sum_{i \in [n]} \frac{A_{ji}}{d_i} = \kappa_j,
\]
it follows from the normalisation that
\[
\Delta_{[n]}^x = \frac{1}{n} \sum_{i \in [n]} x_i \kappa_i - \left[\frac{1}{n} \sum_{i \in [n]} x_i\right]
\left[\frac{1}{n} \sum_{i \in [n]} \kappa_i\right] = \mathrm{Cov}(x_U,\kappa_U).
\]
\vrdg{It follows that $x_U$ and $\kappa_U$ are positively correlated if and only if the GFP holds.}
\end{proof}

\begin{remark}
{\rm Since $\mathrm{Cov}(d_U,\kappa_U) = \Delta^d_{[n]} \geq 0$, we have that $d_U$ and $\kappa_U$ are positively correlated. In \cite{CKN2021} the following statement is made: “While it is not mathematically guaranteed, in practice we expect properties that correlate with $d$ to also correlate with $\kappa$, and hence in such situations we can reasonably expect to observe a generalised friendship paradox.'' In general, this hope is unfounded, because $\mathrm{Cov}(x_U,d_U) \geq 0$ does not imply $\mathrm{Cov}(x_U,\kappa_U) \geq 0$, even though $\mathrm{Cov}(d_U,\kappa_U) \geq 0$.} \hfill$\spadesuit$
\end{remark}

\begin{proposition}
\label{prop:cov}
Suppose that $x_i = d_i f(d_i)$ with $f\colon\,\N \to [0,\infty)$ non-decreasing. Then the GFP holds.
\end{proposition}

\begin{proof}
Write
\[
\sum_{i \in [n]} x_i \kappa_i = \sum_{i \in [n]} d_i f(d_i) \sum_{j \in [n]} \frac{A_{ij}}{d_j}
\]
and
\[
\sum_{i \in [n]} x_i = \sum_{i \in [n]} d_if(d_i),
\]
to get
\[
\mathrm{Cov}(x_U,\kappa_U) = \frac{1}{n} \sum_{i \in [n]} \sum_{j \in [n]} f(d_i) A_{ij} \left(\frac{d_i}{d_j} - 1\right).
\]
Symmetrisation gives
\[
\begin{aligned}
\mathrm{Cov}(x_U,\kappa_U) &= \frac{1}{2n} \sum_{i \in [n]} \sum_{j \in [n]} A_{ij}
\left[f(d_i)\left(\frac{d_i}{d_j} - 1\right) + f(d_j)\left(\frac{d_j}{d_i} - 1\right)\right]\\[0.2cm]
&= \frac{1}{2n} \sum_{i \in [n]} \sum_{j \in [n]} \frac{A_{ij}}{d_id_j}
[d_i f(d_i)- d_j f(d_j)] (d_i-d_j).
\end{aligned}
\]
Since $f$ is non-negative and non-decreasing, the two differences have the same sign, which makes the summand non-negative. Hence $\mathrm{Cov}(x_U,\kappa_U) \geq 0$, and so the GFP holds. Note that, when $f$ is strictly positive, equality holds if and only if all connected components are regular.
\end{proof}

\begin{figure}[htbp]
\centering
\begin{tikzpicture}[scale = 0.5,
  every node/.style={circle, draw, minimum size=10pt, inner sep=0pt, font=\fontsize{2.5}{4}\selectfont}]
    \node (2) at (0,6) {2};
    \node (4) at (0,4) {4};
    \node (1) at (-2,4) {1};
    \node (3) at (2,4) {3};
    \node (5) at (0,2) {5};
    \node (6) at (-2,0) {6};
    \node (7) at (2,0) {7};
    \node (8) at (4,0) {8};
    \node (9) at (4,2) {9};
    \node (10) at (6,0) {10};
    \node (11) at (4,-2) {11};
    \draw (2) -- (4) -- (5) -- (6) -- (7) -- (8) -- (10);
    \draw (1) -- (4);
    \draw (3) -- (4);
    \draw (5) -- (7);
    \draw (9) -- (8);
    \draw (11) -- (8);
\end{tikzpicture}
\vspace{-0.1cm}
\caption{\small A counter example for the TFP: $n=11$, $\Delta_{i,n}^x = 0$ for $i=1,2,3,6,9,10,11$, $\Delta_{i,n}^x = \tfrac14$ for $i=4,8$, $\Delta_{i,n}^x = - \tfrac13$ for $i=5,7$, $\Delta_{[n]}^x = -\frac{1}{66}$. In this example, the degree and the number of triangles do not positively correlate with each other: Vertices 4 and 8 have the largest degrees but are in no triangle.}
\label{fig:TFPcounter}
\end{figure}
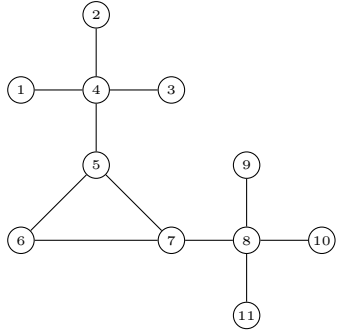

Let $x_i$ be the number of \emph{wedges} that contain vertex $i$. Then $x_i = \tfrac12 d_i(d_i-1)$. Hence, by Proposition~\ref{prop:cov} with
\[
f(d) = \tfrac12(d-1), \quad d \in\N,
\]
the GFP holds for wedges. The number of \emph{triangles} at $i$, however, does not only depend on the edges towards the neighbours of $i$, but also on the edges between the neighbours of $i$. This means that $x_i$ cannot be written as a function of $d_i$. Figure~\ref{fig:TFPcounter} shows a counter example \cite{M2024}.

The mathematical literature on the GFP is scarce \cite{PYNSBN2019}, and much remains unexplored. In the present paper we focus on triangles, i.e., we zoom in on the \emph{Triangle Friendship Paradox} (TFP).

\vspace{-0.4cm}


\subsection{Main \vrdg{results}}
\label{ss.theorems}

\vrdg{Our main results concern the Triangle Friendship Paradox (TFP) for various classes of graphs.

\paragraph{TFP.}
Let $G_n$ be a finite undirected graph with $n$ vertices labelled by $[n] = \{1,\ldots,n\}$. Let $\Delta^T_{i,n}$ denote the \emph{triangle friendship-bias of vertex $i \in [n]$} given by
\begin{equation}
\Delta^T_{i, n} =  \left[\frac{1}{d_i} \sum_{j \in [n]} A_{ij} t_j - t_i\right] 1_{\{d_i \neq 0\}},
\label{eq:tfpbias}
\end{equation}
where $(A_{ij})_{i,j \in [n]}$ is the adjacency matrix, $d_i$ is the degree of vertex $i \in [n]$, and $t_i$ is the number of triangles that contain vertex $i\in [n]$. The \emph{Triangle Friendship Bias} (TFB) is defined as
\begin{equation}
\Delta_{[n]}^T = \frac{1}{n} \sum_{i \in [n]} \Delta_{i,n}^T,
\label{eq:tfbgraph}
\end{equation}
and the TFP holds when $\Delta_{[n]}^T  \ge 0$.}

Section~\ref{sss.star} looks at \emph{deterministic finite graphs}, Section~\ref{sss.sparse} at \emph{sparse random graphs}, and Section~\ref{sss.dense} at \emph{dense random graphs}. Deterministic graphs are harder to handle than random graphs. This is why our results for deterministic graphs are more restrictive than for random graphs. For the latter we can use the notions of \emph{local convergence} (for spare random graphs) and \emph{graphons} (for dense random graphs), which are powerful tools.


\subsubsection{Finite deterministic graphs}
\label{sss.star}

\begin{figure}[htbp]
\centering
\begin{tikzpicture}[scale = 0.4, every node/.style={circle, draw, minimum size=6pt, inner sep=1pt, font=\tiny}
]
\newcommand{\drawKfive}[1]{
\foreach \i in {1,...,5} {
\node (N#1\i) at ({#1 + 2*sin(72*(\i))}, {2*cos(72*(\i))}) {\i};
}
\foreach \i in {1,...,5} {
\foreach \j in {\i,...,5} {
\ifnum\i<\j
\draw[-] (N#1\i) -- (N#1\j);
\fi
}
}
\node (T1#1) at ({#1 + 2*sin(72*(4))}, {4 - 2*cos(72*(4)}) {6};
\node (T2#1) at ({#1 + 2*sin(72*(1))}, {4 - 2*cos(72*(1)}) {7};
\draw[-] (N#15) -- (T1#1) -- (T2#1)--(N#15);
}
\drawKfive{0}
\drawKfive{8}
\drawKfive{16}

\node (A) at ({8 + 4*sin(72*(4))- 2*sin(72*(5))}, {6 - 4*cos(72*(4)}){8};
\node (B) at ({8 + 2*sin(72*(5))}, {6 - 4*cos(72*(4)}){9};
\draw[color = vrr] (A) -- (B);
\draw[color = vrr] (A) -- (T18);
\draw[color = vrr] (T18) -- (B);
\node[fill = blue!20] (T18) at ({8 + 2*sin(72*(4))}, {4 - 2*cos(72*(4)}) {6};

\node (C) at ({16 + 2*sin(72*(0))}, {6 - 4*cos(72*(4)}){8};
\draw[color = vrr] (C) -- (T116);
\draw[color = vrr] (C) -- (T216);
\draw[color = blue] (T116) -- (T216);

\end{tikzpicture}
\caption{\small Two examples where adding a triangle (coloured red) to the first graph reduces the average friendship-bias. The bias drops from $\tfrac{13}{21}$ to $\tfrac{7}{18}$ when the triangle is added at the vertex, as in the second graph. The bias drops from $\tfrac{13}{21}$ to $\tfrac{7}{24}$ when the triangle is added at the edge, as in the third graph.}
\label{fig:addtriangle}

\end{figure}
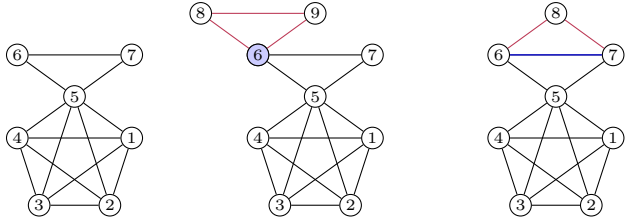

In this section we identify several finite graphs for which the TFP holds. For instance:
\begin{itemize}
\item
The TFP holds for all graphs where every edge is part of exactly one triangle, i.e., the graph is made up of triangles that are joined together at their vertices. Indeed, in that case for each vertex the number of triangles it lies in is exactly half its degree. Hence, $t_i = \tfrac12 d_i$ and, using that the TFP holds for degrees, we get the claim.
\item 
Consider a graph in which the distance between any pair of triangles is at least $3$. In that case the triangle friendship-bias of each vertex is determined by a single triangle. The only vertices with a non-zero bias either fall in a triangle or are at distance $1$ from a triangle, and hence the total bias is given by
\[
n \Delta^T_{[n]} = \sum_{ {i \in [n] \colon} \atop {|i-\triangle| = 0} } \left(\frac{2}{d_i}-1\right)
+ \sum_{ {i \in [n]\colon} \atop {|i-\triangle| = 1} } \frac{1}{d_i},
\]
where $\triangle$ denotes the collection of triangles, and $d_i \geq 2$ in the first sum and $d_i \geq 1$ in the second sum. Suppose that some of the vertices in the triangles are connected to other vertices in other triangles by (possibly multiple) disjoint paths of length $\geq 3$. Then for every $i \in [n]$ such that $|i-\triangle| = 0$ there are precisely $d_i-2$ vertices $j \in [n]$ such that $|j-\triangle| = 1$, each of which has $d_j=2$. Hence
\[
n \Delta^T_{[n]} = \sum_{ {i \in [n] \colon} \atop {|i-\triangle| = 0} } \left(\frac{2}{d_i}-1+\frac{d_i-2}{2}\right)
= \sum_{ {i \in [n] \colon} \atop {|i-\triangle| = 0} } \left(\sqrt{\frac{d_i}{2}} - \sqrt{\frac{2}{d_i}}\,\right)^2 \geq 0,
\]
and so the TFP holds. Strict inequality holds unless $d_i=2$ for all $i \in [n]$ such that $|i-\triangle| = 0$, which corresponds to isolated triangles. Conversely, by attaching many disjoint edges to each of the vertices $i \in [n]$ such that $|i-\triangle| = 0$ or $|i-\triangle| = 1$, we can make $d_i$ arbitrarily large, in which case the total bias is close to $-|\triangle|$ and the TFP fails.
\item
Adding a triangle either at a vertex or at an edge does not always increase the average friendship-bias. Two examples are shown in Figure~\ref{fig:addtriangle}. In both examples the bias is strictly positive both before and after the triangle is added, yet there is a drop.
\end{itemize}

\vspace{-0.3cm}

\vrdg{Our} two main theorems below focus on \emph{partially completed star-graphs}, which turn out to have a tractable structure. \vrdg{We define them before stating our main results.}

\begin{figure}[htbp]
\centering
\begin{tikzpicture}[scale = 0.9, every node/.style={circle, draw, minimum size=9pt, inner sep=0pt, font=\tiny}]

\node (0) at (0,0) {0};

\foreach \i in {1,...,14} {
\node (\i) at ({2*cos(360/14*\i)}, {2*sin(360/14*\i)}) {\i};
\draw (0) -- (\i);
}
\draw[color = vrr, thick] (0) -- (4);
\draw[color = vrr, thick] (0) -- (1);
\draw[color = vrdg, thick] (0) -- (3);
\draw[color = vrdg, thick] (0) -- (2);
\draw[ color = vrdg, thick] (2) -- (3);
\draw[color = vrdg, thick] (0) -- (14);
\draw[color = vrdg, thick] (0) -- (13);
\draw[ color = vrdg, thick] (14) -- (13);
\draw (11) -- (12);
\draw[-] (11) -- (10);
\draw[-] (9) -- (10);
\draw[-] (9) -- (8);
\draw[-] (7) -- (6);
\draw[-] (6) -- (5);
\end{tikzpicture}
\caption{\small A partially completed star-graph $G(2,2,\{2,5\})$ with 15 vertices (including the center), two bands of adjacent triangles of width 2 and 4, two isolated triangles, and two tadpoles. This graph has $\Delta_{[15]}^T = \frac{158}{45}$.}
\label{fig:pcs}
\end{figure}
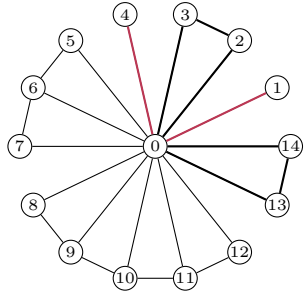

\begin{definition}\label{def:pstar}
{\rm \vrdg{Consider a graph with vertex set $\{0\}\cup [n]$ for $n>1$ and edge set $E = \{\{0,i\}\colon\, i\in [n]\} \cup W$ for some $\varnothing \neq W \subsetneq \{\{i, i+1\}\colon\, i\in [n]\}$, the cyclic ring of edges around $0$ with the convention $n+1=1$ (see Figure~\ref{fig:pcs}). We call an edge $\{0, i\}$ a tadpole if vertex $i\in [n]$ is connected only to vertex $0$. We will refer to this graph as} a partially completed star-graph and denote it by $G(t,\tilde k, \{k_i\}_{i=1}^m)$, where \vrdg{$t$ is the number of} tadpoles, $\tilde k$ \vrdg{is the number of} isolated triangles, and \vrdg{there are} $m$ bands of adjacent triangles with \vrdg{$k_i$ denoting number of triangles in the $i$th band}. We will refer to vertex 0 as the centre, and \vrdg{denote total number of triangles by} $k = \tilde k + \sum_{i=1}^m k_i$. Note that $n = 2\tilde k + \sum_{i=1}^m (k_i+1)+t$.}\hfill$\spadesuit$
\end{definition}
We agree that the triangle with vertex set $\{0,1,2\}$ and edge set $\{\{0,1\}, \{0,2\}, \{1,2\}\}$ is a partially completed star-graph, although it does not strictly satisfy Definition~\ref{def:pstar}.

\vrdg{We are now ready to state our first main result on partially completed star-graphs.}

\begin{theorem}
\label{thm:stargraphs}
\vrdg{The TFP holds for partially completed star-graphs. For $G(t,\tilde k, \{k_i\}_{i=1}^m)$ the TFB is given by}
\begin{equation}
\label{PStarTbiasmain}
\Delta^T_{[n+1]} = \frac{1}{n+1}\left(\C + \sum_{i =1}^m \AT_k(k_i) + \tilde k\, \IT_k + t\, \TP_k\right),
\end{equation}
where
\vrdg{
\begin{align}
&\C = \tfrac{2k}{n}-k &\mbox{ is the TFB of the center $0$,}\hspace{3.6cm}\nonumber \\
&\AT_k(k_i) = \tfrac{k(k_i+2)}{3} - \tfrac{2k_i}{3}
&\mbox{ is the sum of the TFB of the non-center vertices}\nonumber \\
&&\mbox{in the band of $k_i$ adjacent triangles,}\hspace{1.5cm}\nonumber\\
&\IT_k = k-1 &\mbox{ is the sum of the TFB of the non-center vertices}\nonumber\\
& &\mbox{in each of the $\tilde{k}$ isolated triangles and}\hspace{1.5cm}\nonumber \\
&\TP_k = k &\mbox{is the TFB of the non-center vertex of a tadpole.}
\label{eq:PCStr}
\end{align}}
\end{theorem}

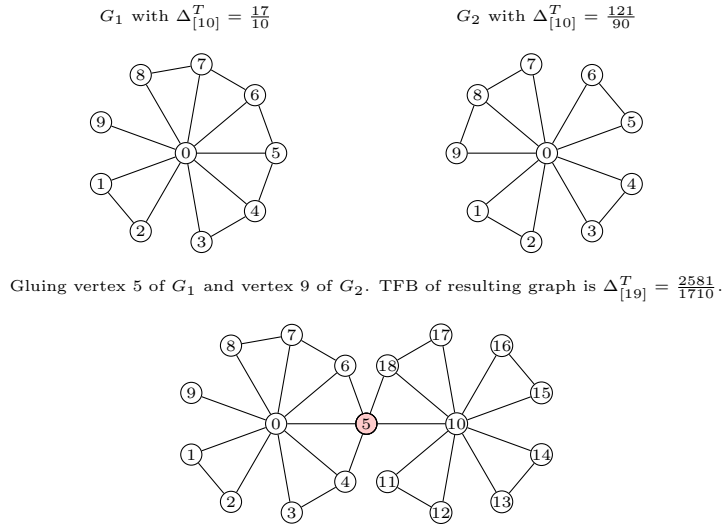
\begin{figure}[htbp]
\centering
\begin{tikzpicture}[scale = 0.6, every node/.style={circle, draw, minimum size=8pt, inner sep=0pt, font=\tiny}]

\node[draw =none] at (0,3) {$G_1$ with $\Delta^T_{[10]} = \frac{17}{10}$};
\node (A0) at (0,0) {0};

\foreach \i in {1,...,9} {
  \node (A\i) at ({2*cos(360/9*(\i+4))}, {2*sin(360/9*(\i+4))}) {\i};
  \draw (A0) -- (A\i);
}
\draw[-] (A1) -- (A2);
\draw[-] (A3) -- (A4) -- (A5) -- (A6) -- (A7) -- (A8);

\node[draw =none] at (8,3) {$G_2$ with $\Delta^T_{[10]} = \frac{121}{90}$};
\node (B0) at (8,0) {0};

\foreach \i in {1,...,9} {
\node (B\i) at ({8+2*cos(360/9*\i+ 180)}, {2*sin(360/9*\i + 180)}) {\i};
\draw (B0) -- (B\i);
}
\draw[-] (B1) -- (B2);
\draw[-] (B3) -- (B4);
\draw[-] (B5) -- (B6);
\draw[-] (B7) -- (B8) -- (B9);

\node[draw = none] at (4, -3) {
Gluing vertex $5$ of $G_1$ and vertex $9$ of $G_2$. TFB of resulting graph is $\Delta^T_{[19]} = \frac{2581}{1710}$.};
\node[fill = red!20] (C0) at (4, -6) {};

\node (C10) at (2,-6) {0};

\foreach \i in {1,...,9} {
\node (C1\i) at ({2+ 2*cos(360/9*(\i+4))}, {2*sin(360/9*(\i+4)) - 6}) {\i};
\draw (C10) -- (C1\i);
}
\draw[-] (C11) -- (C12);
\draw[-] (C13) -- (C14) -- (C15) -- (C16) -- (C17) -- (C18);

\node (C20) at (6,-6) {10};

\foreach \i in {11,...,18} {
\node (C2\i) at ({6+2*cos(360/9*(\i-10)+ 180)}, {2*sin(360/9*(\i-10) + 180)-6}) {\i};
\draw (C20) -- (C2\i);
}

\node (C29) at (4, -6) {};
\draw (C20) -- (C29);
\draw[-] (C211) -- (C212);
\draw[-] (C213) -- (C214);
\draw[-] (C215) -- (C216);
\draw[-] (C217) -- (C218) -- (C29);
\end{tikzpicture}
\vspace{-1.5cm}
\caption{\small Illustration of Theorem \ref{thm:stargraphsglued}. Two partially completed star-graphs $G_1$ and $G_2$ of size 10 each are glued together at a vertex that is not the center.}
\label{fig:StarConcatenate_plot}
\end{figure}

\vrdg{Our next result shows that} the TFP is preserved when we glue two partially completed star-graphs \vrdg{together by merging two non-center vertices from each graph (with edge set being the union of the two edge sets). See Figure \ref{fig:StarConcatenate_plot} for an illustration.}

\begin{theorem}
\label{thm:stargraphsglued}
Let $G_{n_1}$ \vrdg{with vertex set $\{0\}\cup [n_1]$} and $G_{n_2}$ \vrdg{with vertex set $\{0\}\cup [n_2]$} be two partially completed star-graphs \vrdg{such that each of them has at least one triangle.} Suppose that $G_{n_1}$ and $G_{n_2}$ are glued together at a vertex that is not the center. Then the glued graph again satisfies the TFP.
\end{theorem}

\vspace{-1.5cm}

It is natural to wonder whether simultaneous gluing of \emph{more than two} partially completed star-graphs also preserves the TFP. The answer is no. Take four \vrdg{partially completed star-graphs each of which is a} triangle with one edge attached to one of the vertices forming a tadpole. Glue the \vrdg{four endpoints of the tadpole} edges together into a single vertex. Each individual \vrdg{partially completed star-graph} has $\Delta^T_{[4]} = \tfrac{1}{6}>0$, but the glued graph, which is not a partially completed star-graph, has $\Delta^T_{[13]} = -\tfrac{1}{39} < 0$.


\subsubsection{Sparse random graphs}
\label{sss.sparse}

\vrdg{In this section, we consider a sequence of sparse random graphs $(G_n)_{n\in\N}$ and study the limiting value of $\Delta_{[n]}^T$ under local convergence. We follow the notation and definitions from \cite[Chapter 2]{vdH2024}. Let $\mathcal{G}$ be the space of all rooted (possibly infinite) undirected graphs modulo isomorphism. The following definition is taken from \cite[Page 56, Definition 2.11]{vdH2024}.

\begin{definition}
\label{deflc}
{\rm Let $(G_{n})_{n\in\N}$ be a sequence of finite random graphs in $\mathcal{G}$. We say that $G_{n}$ converges \emph{locally in probability} to $(G,o)\in \mathcal{G}$, with (a possibly random) law $\nu$ when, for every bounded and continuous function $h\colon\,\mathcal{G}\to\R$,
\begin{align*}
\mathbb{E}\big[h(G_{n} ,U_{n}) \mid G_{n}\big] \quad \mbox{converges in probability to }
\quad \mathbb{E}_{\nu}\big[h(G,o)\big] \quad \text{as $n\to \infty$},
\end{align*}
where $\mathbb{E}$ is expectation with respect to the random vertex $U_{n}$ and the random graph $G_{n}$.
}\hfill$\spadesuit$
\end{definition}
}

\vrdg{Our first result states that if the graphs converge locally in probability, then the corresponding TFBs converge in distribution.}

\begin{theorem}
\label{thm:sparse}
If $(G_n)_{n\in\N}$ converges locally in probability to $(G,o)$, then $\Delta^T_{[n]}$ converges in distribution to \vrdg{a random variable} $\Delta^T$, defined to be the triangle friendship-bias of $o$ in $G$.
\end{theorem}

Many sparse random graphs are \emph{locally tree-like}, and therefore locally converge to a \emph{rooted random tree} (see \cite[Chapters 3--5]{vdH2024}). On trees, all triangle counts are zero, and so $\Delta^T= 0$ a.s., so that $\lim_{n\to\infty} \Delta^T_{[n]} = 0$. Below we look at two examples for which we can identify how $\Delta^T_{[n]}$ scales as $n\to\infty$.

\vrdg{We now define two sequences of random graphs namely, the Erd\H{o}s-R\'enyi random graph and the configuration model. We will assume without loss of generality that they are defined on a common probability space $(\Omega, \mathcal{F}, \mathbb{P})$, and $\mathbb{E}$ will denote expectation with respect to $\mathbb{P}$.}


\paragraph{Sparse Erd\H{o}s-R\'enyi Random Graph.}

\begin{figure}[htbp]
\begin{center}
\vspace{-1cm}
\includegraphics[width=.3\textwidth]{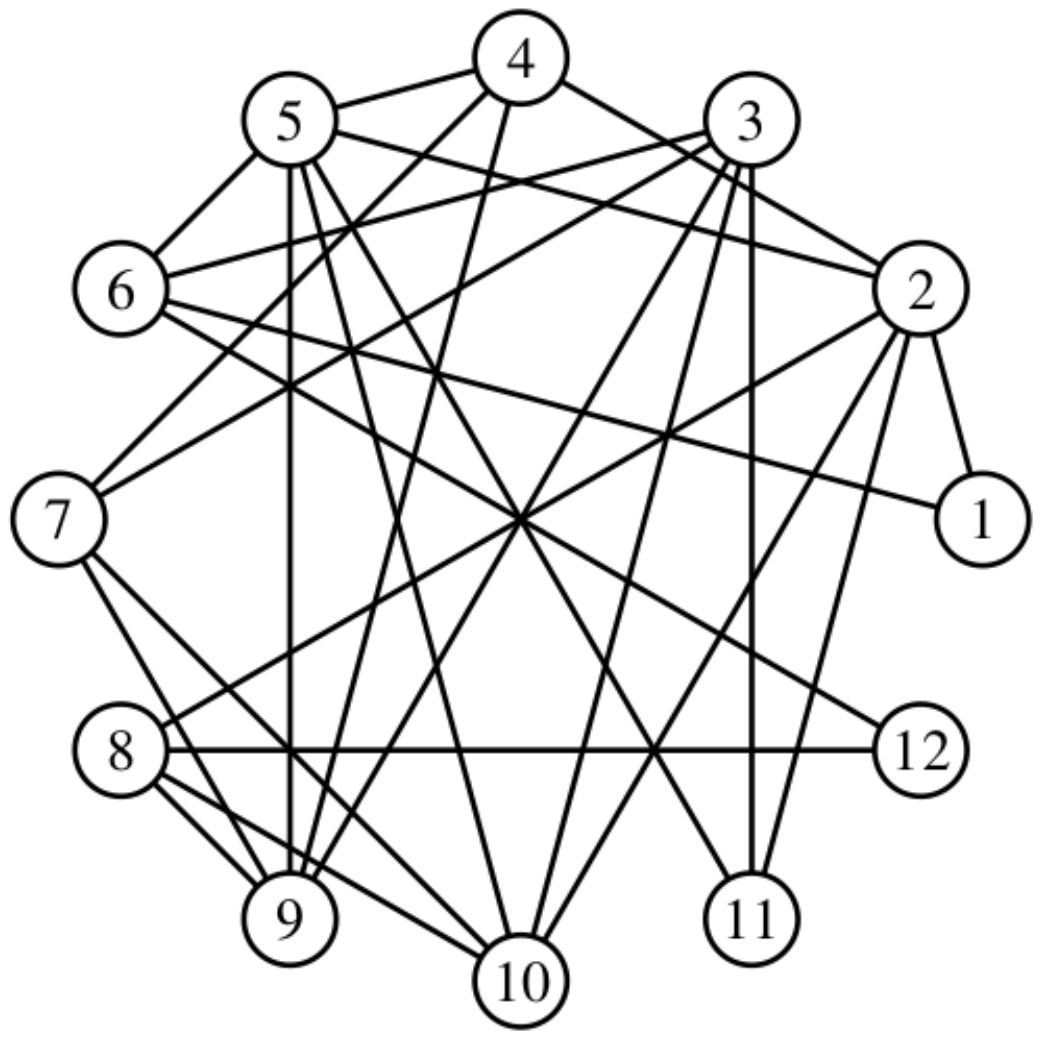}
\end{center}
\vspace{-1.5cm}
\caption{\small A realisation of an ERRG$(12, \frac{1}{3})$ with $n = 12$ and $\lambda = 4$.}
\label{fig:ERRG}
\end{figure}

The ERRG$(n, \frac{\lambda}{n})$ with $n$ vertices and edge retention probability $\lambda/n \in (0,1]$, with $\lambda \in (0,\infty)$ and $n \geq \lambda$, is the sparse random graph obtained from the complete graph $K_n$ by independently deleting edges with probability $1-(\lambda/n)$. See Figure \ref{fig:ERRG} for an illustration.

\begin{theorem}
\label{thm:ER}
For the sparse ERRG$(n, \frac{\lambda}{n})$ the following hold:
\begin{itemize}
\item[(a)]
For every $n \geq 3$ and $\lambda \in (0,\infty)$,
\begin{equation}
\mathbb{E}\left[\Delta_{[n]}^T\right] = \tfrac{\lambda}{n}\left(1-\tfrac{\lambda}{n}\right)\,\left[\tfrac{\lambda}{n}(n-2)-1\right]
+ \left[\tfrac{\lambda}{n}-(\tfrac{\lambda}{n})^3\tfrac{1}{2}(n-2)(n-3)\right]\,\left(1-\tfrac{\lambda}{n}\right)^{n-1}.
\label{eq:meanERRG}
\end{equation}
\item[(b)]
For every $\lambda \in (0,\infty)$,
\begin{equation}
\lim_{n\to\infty} \mathbb{E}\left[n\Delta_{[n]}^T\right] = \zeta(\lambda)
\label{eq:scaledmean}
\end{equation}
with $\zeta(\lambda) = \lambda^2 - \lambda + [\lambda-\tfrac12\lambda^3]\,\eee^{-\lambda}$.
\item[(c)]
For every $\lambda \in (0,\infty)$, the weak law of large numbers fails because
\begin{equation}
\liminf_{n\to\infty} \mathbb{P}\left(n\Delta_{[n]}^T = 0\right) > 0.
\label{eq:wlf}
\end{equation}
\end{itemize}
\end{theorem}

\begin{figure}[htbp]
\vspace{-2cm}
\centering
\includegraphics[width =0.4\textwidth]{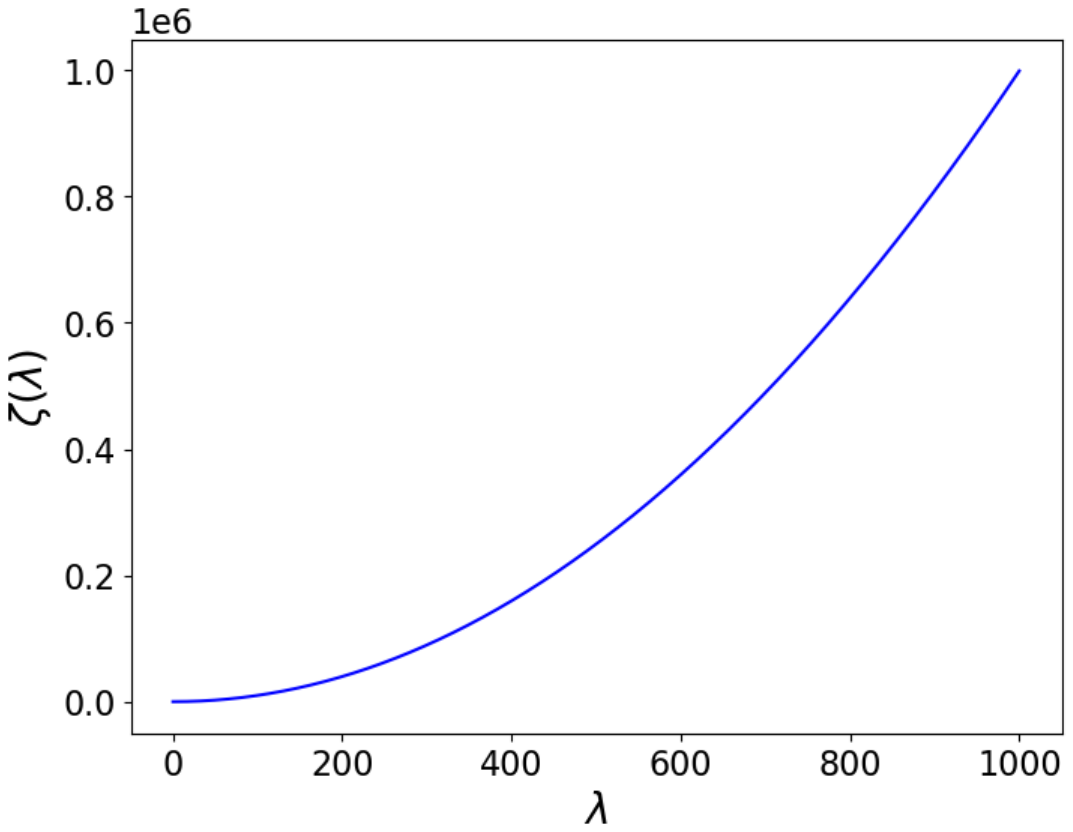}
\vspace{-2cm}
\caption{\small Numerical plot of $\lambda \mapsto \zeta(\lambda)$ for $10^{-3} \leq \lambda \leq 10^3$.}
\label{fig:zetalambda_plot}
\end{figure}

\noindent
Figure \ref{fig:zetalambda_plot} shows a numerical plot $\lambda \mapsto \zeta(\lambda)$, which is strictly increasing (a proof of this property is given in Remark~\ref{rem:increasing}). Thus, \vrdg{the expected value of TFB is positive} for large enough sparse Erd\H{o}s-R\'enyi random graphs, and \vrdg{their limit} is a strictly increasing function of the edge density. Note that $\zeta(\lambda) \sim \tfrac12 \lambda^4$ as $\lambda \downarrow 0$ and $\zeta(\lambda) \sim \lambda^2$ as $\lambda\to\infty$.

We will see that the number of triangles is finite and has distribution $\mathrm{POISSON}(\lambda^3/6)$. It is possible to show that $\Delta_{[n]}^T/\mathbb{E}[\Delta_{[n]}^T]$ converges in law to a random variable $Z$, and to identify the probability distribution of $Z$. The latter can be expressed in terms of independent Poisson random variables. \vrdg{We do not provide a proof here, as our focus is on the TFP itself.}


\paragraph{Configuration Model.}

\begin{figure}[htbp]
\vspace{0.1cm}
\centering
\begin{tikzpicture}[every loop/.style={}]
\foreach \ang [count=\n from 1] in {60,0,...,-240}
\node (v\n) [cnd] at (\ang:1.5cm) {$\n$};
\draw [-Bar] (v1.240) -- +(240:0.45);
\draw [-Bar] (v2.-180) -- +(180:0.45);
\draw [-Bar] (v2.-225) -- +(-225:0.45);
\draw [-Bar] (v2.-135) -- +(-135:0.45);
\draw [-Bar] (v4.80) -- +(80:0.45);
\draw [-Bar] (v4.40) -- +(40:0.45);
\draw [-Bar] (v4.60) -- +(60:0.45);
\draw [-Bar] (v5.30) -- +(30:0.45);
\draw [-Bar] (v5.-30) -- +(-30:0.45);
\draw [-Bar] (v6.-20) -- +(-20:0.45);
\draw [-Bar] (v6.-40) -- +(-40:0.45);
\draw [-Bar] (v6.-60) -- +(-60:0.45);
\draw [-Bar] (v6.-80) -- +(-80:0.45);
\draw [-Bar] (v3.120) -- +(120:0.45);
\begin{scope}[xshift=5cm]
\foreach \ang [count=\n from 1] in {60,0,...,-240}
\node (v\n) [cnd] at (\ang:1.5cm) {$\n$};
\draw [-] (v1) -- (v6) -- (v5) -- (v3);
\draw [-] (v2) -- (v4);
\draw (v6) to[bend right = 15] (v2);
\draw (v6) to[bend left = 15] (v2);
\draw (v4) to[out=180, in=145, looseness=16] (v4);
\end{scope}
\end{tikzpicture}
\caption{\small Pairing algorithm for the Configuration Model with $n=6$ and degree sequence $(1,3,1,3,2,4)$. The vertices are labelled clockwise from the right-top. Vertex $i$ is assigned $d_i$ half-edges. The half-edges are paired uniformly at random to become edges.}
\label{fig:pairing}
\end{figure}
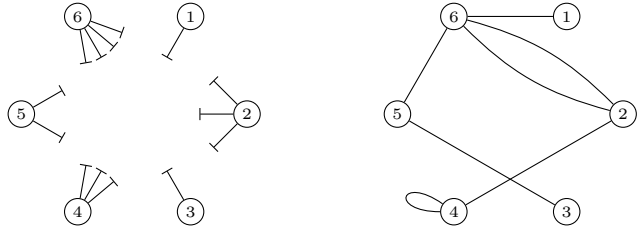

The $\cm$ with $n$ vertices is a graph with a \emph{prescribed degree sequence} $\vrdg{\mathbf{d} = }(d_i)_{i \in [n]} \in \N^n$ that can be generated via a pairing algorithm in which half-edges are matched randomly to obtain edges. The matching leads to a \emph{random multi-graph} with the property that the degree of vertex $i$ is $d_i$ for every $i \in [n]$. The total number of edges is $\tfrac12\sum_{i \in [n]} d_i$, and so the sum of the prescribed degrees must be even. See \cite[Chapter 7]{vdH2017} for precise definitions and Figure~\ref{fig:pairing} for an illustration.

\vrdg{\paragraph{Convention:}The random graph $\cm$ allows for self-loops and multiple edges. We adopt the \emph{convention} that the expression in  \eqref{eq:tfpbias} is the TFB of $\cm$ at vertex $i\in [n]$. To be precise, the multiple edges add to the triangle count and there is an additional term $A_{ii}t_i$ in \eqref{eq:tfpbias} (since $A_{ii}$ may be non-zero). To simplify computations and scaling involved, we will work with \eqref{eq:tfpbias} and the convention that the TFB for $\cm$ is given by \eqref{eq:tfbgraph}.}

For $n \in \N$, define
\[
m_k = m_k(n) = \sum_{i \in [n]} d_i^k, \qquad k \in \N.
\]
Suppose that
\begin{equation}
\label{CMscal}
c_k = \lim_{n\to\infty} n^{-1}m_k(n) \text{ exists and is strictly positive and finite for } 1 \leq k \leq 7.
\end{equation}

\begin{theorem}
\label{thm:CM}
For the sparse $\cm$ the following hold:
\begin{itemize}
\item[(a)]
If $m_1>7$, then
\[
\mathbb{E}[\Delta_{[n]}^T]
 = \frac{1}{2(m_1-1)(m_1-3)(m_1-5)(m_1-7)}
\vrdg{\left(\frac{A_1+A_2+A_3}{n} + (B_1+B_2+B_3)\right)},
\]
with
\begin{align*}
\vrdg{
A_1} &\vrdg{= 6m_1m_3 - 18 m_1m_4 + 14 m_1m_5 - 2 m_1m_6 + 4 m_2m_3 - 8 m_2m_4 + 4 m_3^2}, \\
\vrdg{A_2} &\vrdg{= - 14 m_1m_2m_3 + 3 m_1m_2m_4 -9 m_1^2m_2 + 18 m_1^2m_3
- 5 m_1^2m_4 +3 m_1m_2^2 + 4 m_2^3}, \\
A_3 &= - 7 m_1^3m_2 - m_1m_2^3 + 5 m_1^2m_2^2 + 3 m_1^4,\\
\intertext{and}
B_1 &= \vrdg{2 m_7 - 10 m_6 + 14 m_5 - 6 m_4},\\
B_2 &= \vrdg{6 m_1m_3 - 8 m_1m_4 + 2 m_1m_5 -13 m_2m_3 + 11 m_2m_4}, \\
&\quad\; \vrdg{- 2 m_2m_5 - m_3m_4 + 3 m_2^2 + 2 m_3^2,}\\
B_3 &=  - 2m_1m_2m_3 - 3 m_1^2m_2 + m_1^2m_3 + 6 m_1m_2^2 +m_2^2m_3 - 3 m_2^3.
\end{align*}
\item[(b)]
Subject to \eqref{CMscal},
\[
\lim_{n\to\infty} \mathbb{E}[n\Delta_{[n]}^T] = \zeta(c_1,c_2,c_3)
\]
with
\[
\zeta(c_1, c_2, c_3) = \frac{(c_2-c_1)^2}{2c_1^4} \big[(c_3-c_1c_2) - 3(c_2-c_1^2)\big] \geq 0.
\]
Equality holds if and only if $c_2 = c_1^2$, i.e., for the special case of the $d$-regular random graph with $c_k=d^k$, $k \in \N$.
\item[(c)]
Subject to \eqref{CMscal}, with the exception of the $d$-regular random graph, the weak law of large numbers fails because
\[
\liminf_{n\to\infty} \mathbb{P}(n\Delta_{[n]}^T = 0) > 0.
\]
\end{itemize}
\end{theorem}

\noindent
The expression for $\mathbb{E}[\Delta_{[n]}^T]$ for finite $n$ is somewhat baroque, but it simplifies considerably in the limit as $n\to\infty$. In all cases where $\zeta(c_1,c_2,c_3)$ is strictly positive, the TFP holds \emph{on average} for large enough sparse $\cm$s.

We will see that, in the limit as $n\to\infty$, the number of triangles is finite and has distribution $\mathrm{POISSON}([(c_2-c_1)/c_1]^3/6)$. All triangles are isolated, so that their contributions to the total \vrdg{triangle} friendship-bias can be summed up. For the $d$-regular random graph each of these contributions is $0$. Indeed, each triangle has 3 vertices, each of which has 2 neighbours inside the triangle and $d-2$ neighbours outside the triangle. Each of these vertices has triangle friendship-bias $\frac{2}{d}-1 = - \frac{d-2}{d}$, while each of its $d-2$ neighbours outside the triangle has triangle friendship-bias $\frac{1}{d}$.

It is again possible to show that $\Delta_{[n]}^T/\mathbb{E}[\Delta_{[n]}^T]$ converges in law to a random variable $Z$, and to identify the probability distribution of $Z$. However, the latter is again rather cumbersome to spell out, and therefore we refrain from doing so.


\subsubsection{Dense random graphs}
\label{sss.dense}

\vrdg{In this section we consider dense random graph sequences. It is well known that graphons arise as limits of dense graphs in the cut-metric, and that finite graphs can be represented as graphons as well. We describe the preliminaries required for our results below. See \cite{L2012} for a detailed and comprehensive treatment of dense graph theory.}

\vrdg{We refer to the set}
\[
\W = \big\{\kappa\colon [0,1]^2 \to [0,1] \colon\,\kappa \text{ is measurable and symmetric}\big\}
\]
as the set of \emph{graphons}. \vrdg{We may think of $\kappa(\ddd x, \ddd y)$ as the fraction of edges with endpoints in $\ddd x$ and $\ddd y$. We may also construct a canonical finite random graph sequence corresponding to $\kappa$ as follows.} Let \vrdg{$G(n, \kappa)$} be the random graph with vertex set $V_n = [n] = \{1,\ldots,n\}$ and edge set $E_n$ obtained by independently connecting vertices $i = \lceil xn \rceil$ and $j = \lceil yn \rceil$ with probability $\kappa(x,y)$. This is also known as the \emph{inhomogeneous Erd\H{o}s-R\'enyi random graph} with \emph{reference graphon} $\kappa$. \vrdg{It can be shown that $\{G(n, \kappa)\}_{n\in\N}$ converges to $\kappa$ in the cut-metric.

Conversely, for any finite graph $G_n$, one may associate an \emph{empirical} graphon $\kappa^{G_n}\colon\,[0,1]^2 \to [0,1]$ given by
\[
\kappa^{G_n}(x,y) = \left\{\begin{array}{ll}
1  & \mbox{if there is an edge between} \lceil xn \rceil \text{ and } \lceil yn \rceil,\\
0  & \mbox{otherwise.}
\end{array}
\right.
\]}

Since graphs are invariant under relabelling of their vertices, it is natural to replace $\W$ by the set of \emph{equivalence classes} $\tilde\W$ of graphons, in which two graphons $\kappa_1,\kappa_2 \in \W$ are said to be equivalent when there exists a Lebesgue-measure-preserving transformation $\sigma$ of $[0,1]$ such that $\kappa_1(\sigma x,\sigma y) = \kappa_2(x,y)$ for all $x,y \in [0,1]$. The advantage of working with $\tilde\W$ over $\W$ is that the former is a compact set \vrdg{under the \emph{cut-metric}}. \vrdg{Our main results are on the almost sure limiting value of TFB for the finite random graph sequence $\{G(n, \kappa)\}_{n\in\N}$ arising from $\kappa \in \mathcal{W}$. We note that the TFB for such sequences are invariant under relabelling of the vertices.}

\begin{theorem}
\label{thm:dense}
\vrdg{Assume $\kappa \in \W$ is continuous and
\[
\int_{[0,1]} \ddd y\,\kappa(x,y) > 0 \quad \forall\,x \in [0,1].
\]
For all $n\ge 1$, let $\Delta^T_{[n]}$ be the TFB for the graph $G(n, \kappa)$. Then, }
$$\lim_{n\to\infty} n^{-2} \Delta^T_{[n]} = \chi^T$$ $\PP$-a.s.\ with
\begin{equation}
\label{eq:chi}
\chi^T = \int_{[0,1]} \ddd x\, \left[\frac{1}{\mathcal{D}(x)} \int_{[0,1]} \ddd y\,\kappa(x,y) \mathcal {T}(y) - \mathcal{T}(x)\right],
\end{equation}
where
\begin{equation}
\mathcal{D}(x) = \int_{[0,1]} \ddd y\,\kappa(x,y), \quad
\mathcal{T}(x) = \tfrac12 \int_{[0,1]^2} \ddd y\,\ddd z\, \kappa(x,y)\kappa(y,z)\kappa(z,x),
\label{eq:dt}
\end{equation}
denote the \vrdg{degree} density, respectively, the triangle density at $x \in [0,1]$.
\end{theorem}

\vrdg{Our next result states that for rank-1 reference graphons the TFP holds, while there exist non-rank-1 reference graphons for which it fails.}

\begin{theorem}
\label{thm:densefail}
\vrdg{Assume that $\kappa \in \mathcal{W}$ as in Theorem \ref{thm:dense}.}
\begin{itemize}
\vrdg{\item[(a)]
If the reference graphon is of rank 1, i.e., $\kappa(x,y) = \nu(x)\nu(y)$ for some $\nu\colon\,[0,1] \to [0,1]$ measurable, strictly positive and continuous, then $\chi^T \geq 0$ with equality if and only if $\nu$ is constant.}
\item[(b)]
\vrdg{There exist} two-block graphons for which $\chi^T<0$.
\end{itemize}
\end{theorem}
\noindent

\vrdg{After the proof of Theorem~\ref{thm:densefail}, in Remark~\ref{rem:sbm} below, we present several examples of stochastic block models for which $\chi^T$ is positive and negative.} We next discuss two specific examples of dense random graph sequences, namely, ERRG$(n, p)$ and $\cm$, and compute their expected TFB.


\paragraph{Dense Erd\H{o}s-R\'enyi Random Graph.}

\vrdg{The ERRG$(n,p)$ is a graph with $n$ vertices and edge retention probability $p \in (0,1)$. It may be viewed as the random graph obtained from the complete graph $K_n$ by independently deleting edges with probability $1-p$. It is standard to note that this corresponds to a dense graph sequence arising from $\kappa(x,y) = p$ for all $x,y\in [0,1]$. From \eqref{eq:chi} it is easy to see that for the sequence $\{\mathrm{ERRG}(n,p)\}_{n\in\N}$ we have $\chi^T = 0$.} From \eqref{eq:errgtfb} in the proof of Theorem~\ref{thm:ER}(a), we have that
\begin{equation}
\label{eq:ERRG}
\mathbb{E}[\Delta_{[n]}^T] = p(1-p)\,\big[p(n-2)-1\big] + \big[p-\tfrac{1}{2}(n-2)(n-3)p^3\big]\,(1-p)^{n-1}.
\end{equation}
In particular, $\lim_{n\to\infty} \mathbb{E}[n^{-1}\Delta_{[n]}^T] = p^2(1-p)$, which shows that the expectation of the average triangle friendship-bias is proportional to $n$ rather than $n^2$. \vrdg{Due to dense nature of the ERRG$(n,p)$, it is natural to expect that the fluctuations of the triangle friendship-biases of the vertices become small on the scale of their expectation. Indeed,} a long but straightforward computation (available from the authors upon request) shows that $\lim_{n\to\infty} \mathrm{Var}[n^{-1}\Delta_{[n]}^T] = 0$. Consequently, $n^{-1}\Delta_{[n]}^T$ satisfies the weak law of large numbers, i.e., the average triangle friendship-bias is eventually positive with a probability tending to $1$ as $n\to\infty$. We believe that also the strong law of large numbers holds, but a proof would require more delicate concentration estimates.


\paragraph{Dense Configuration Model.}

In the $\cm$, we can choose the degrees to be of order $n$ to get a dense random graph. The formula in Theorem 1.10(a) is valid for any choice of degrees as long as $n$ is finite. If we assume that the degrees are chosen such that
\[
\lim_{n \to \infty}\frac{m_k(n)}{n^{k+1}} = c_k^* \text{ exists and is strictly positive and finite for } 1 \leq k \leq 7,
\]
then we get
\[
\lim_{n \to \infty} \mathbb{E}\big[n^{-2}\Delta_{[n]}^T\big]
= \frac{(c_2^*)^2(c_3^*- c_1^*c_2^*)}{2(c_1^*)^4} \geq 0,
\]
with equality if and only if $c_3^* = c_1^* c_2^*$. For the \emph{regular} $\cm$ with $d_i = d^* n$, $i \in [n]$, we have $c_k^* = (d^*)^k$, and the limit equals 0. When $c_3^* > c_1^* c_2^*$, the expectation of the average triangle friendship-bias is proportional to $n^2$ rather than $n$. Apparently, in that case the behaviour of the dense $\cm$ is similar to that of the dense inhomogeneous \vrdg{ERRG sequence, $\{G(n,\kappa)\}_{n\in \N}$} with a non-constant reference graphon \vrdg{$\kappa$}, and dissimilar to that of the dense homogeneous \vrdg{sequence $\{\mathrm{ERRG}(n,p)\}_{n\in \N}$} with a constant reference graphon \vrdg{$\kappa \equiv p$}.

In principle, $\mathrm{Var}[n^{-2}\Delta_{[n]}^T]$ can be computed to verify whether or not it tends to zero as $n \to \infty$. However, this would require a lengthy computation, and we refrain from doing so.


\subsection{Discussion}
\label{ss.disc}

Our theorems in Section \ref{ss.theorems} describe the average triangle friendship-bias for three classes of graphs: (1) deterministic; (2) random and sparse; (3) random and dense. They provide conditions under which the TFP holds and conditions under which the TFP fails. Since the friendship paradox always holds for degrees and for wedges, it is interesting to see in what ways triangles behave differently. Our theorems offer a \emph{partial classification} and \emph{quantification}.

Deterministic graphs are much harder to analyse than random graphs. Different scenarios appear to be possible. Already the analysis of partially completed star-graphs and their gluing is technically involved, and reveals the intrinsic difficulties associated with handling triangles in this context. For random graphs there seems to be a tendency for the TFP to hold for homogenous and weakly inhomogeneous random graphs and to fail for strongly inhomogeneous random graphs, but there are exceptions.

Challenges for the future are: (1) Show that all graphs consisting of an arbitrary concatenation of triangles satisfy the TFP. (2) Analyse what happens for sparse random graphs with heavy-tailed degree distributions, in particular, preferential attachment random graphs. (3) Find the necessary and sufficient conditions on the reference graphon under which the TFP holds for dense random graphs. Quantify how common it is for the TFP to hold.

\vrdg{\paragraph{Outline of the remainder of the paper.} In Section~\ref{s.proof1}, we prove our two main results for partially completed star-graphs (Theorem~\ref{thm:stargraphs} and Theorem~\ref{thm:stargraphsglued}). In Section~\ref{s.proof2}, we prove that the TFB converges weakly when the graph converges locally in probability (Theorem~\ref{thm:sparse}). Using this, we prove our two main results concerning the TFP for the sparse Erd\H{o}s-R\'enyi random graph (Theorem~\ref{thm:ER}) and the sparse Configuration Model (Theorem~\ref{thm:CM}). In Section~\ref{s.proof3}, we first derive the limiting expression for the TFB for a general dense random graph that arises from a given reference graphon (Theorem~\ref{thm:dense}), and afterwards provide examples for which the limiting TFB is positive and examples for which it is negative (Theorem~\ref{thm:densefail}).}


\section{Proofs for partially completed star-graphs}
\label{s.proof1}

\vrdg{In this section we prove Theorems \ref{thm:stargraphs} and \ref{thm:stargraphsglued}. The proof of Theorem~\ref{thm:stargraphs} is a standard computation of TFB of each vertex according to its type (center, end of a tadpole, endpoint of a triangle). The proof of Theorem~\ref{thm:stargraphsglued} is more involved, and requires that we distinguish several cases and derive a few preliminary lemmas. Lemma~\ref{lem:list1.5} identifies partially completed star-graphs on $\{0\}\cup [n]$ for which the sum of the TFB of each vertex, given by $(n+1)\Delta^T_{[n+1]}$, is less than $\frac{3}{2}$. This result is needed along the way and is of independent interest.}


\subsection{Proof of Theorem \ref{thm:stargraphs}}

\begin{proof}
\vrdg{Let $G(t, \tilde{k}, \{k_i\}_{i=1}^m)$ be a partially completed star-graph on $\{0\}\cup [n]$, and let the total number of triangles be denoted by $k = \tilde k + \sum_{i=1}^m k_i$.} The center $0$ is attached to all the other $n$ non-center vertices. The sum of the number of triangles the $n$ non-center vertices are part of equals $2k$, since each triangle is counted twice, once for each non-center endpoint. \vrdg{Therefore the TFB of vertex $0$ denoted by $C$ is given by}
\begin{equation}
\label{eq:C}
\C = \tfrac{2k}{n} - k.
\end{equation}
\vrdg{For computing the TFB of other vertices, the $k$ triangles that vertex $0$ is attached to will play a central role and we will include the symbol $k$ in all the notation going forward. Fix an $1\le i\le m$, and} a band of $k_i$ adjacent triangles in the graph. \vrdg{We classify} each of the \vrdg{non-center} $k_i+1$ vertices as the two vertices at the end points, the penultimate vertices and the remaining vertices. \vrdg{We compute TFB for each such vertex and sum the outcome to obtain the total TFB for the $k_i + 1$ non-center vertices as}
\begin{equation}
\begin{aligned}
\label{eq:AT}
\AT_k(k_i)
&= 2\left(\tfrac{k+2}{2}-1\right) + \left(\mathbbm{1}_{\{k_i= 2\}} + 2\mathbbm{1}_{\{k_i > 2\}}\right)
\left(\tfrac{k+2 + \mathbbm{1}_{\{k_i>2\}}}{3}-2\right) + \left(k_i - 3\right) \mathbbm{1}_{\{k_i \ge 4\}}
\left(\tfrac{k+4}{3}-2\right)\\
&= \tfrac{k(k_i+2)}{3} - \tfrac{2k_i}{3}.
\end{aligned}
\end{equation}
For \vrdg{each of the $\tilde{k}$} isolated triangles, the sum of the TFBs of both the non-center vertices is
\begin{equation}
\label{eq:IT}
\IT_k = 2 \left(\tfrac{k+1}{2}-1\right) = k-1.
\end{equation}
The TFB of \vrdg{a non-center vertex of} a  tadpole is
\begin{align}
\label{eq:TP}
\TP_k = \tfrac{k}{1}-0 = k.
\end{align}
\vrdg{Therefore, from the definition of the TFB,
\begin{equation}
\Delta^T_\vrdg{{n+1}} = \frac{1}{n+1}\left(\C + \sum_{i =1}^m \AT_k(k_i) + \tilde k\, \IT_k + t\, \TP_k\right).
\label{eq:eqnrepeat}
\end{equation}
By \eqref{eq:C}--\eqref{eq:TP}, \eqref{eq:PCStr} is now immediate.}

\vrdg{Next, we show that the TFP holds.} Since we have assumed that $k\geq 1$, we know that $\AT_k(k_i)$ \vrdg{for all $1\le i\le m$}, $\IT_k$ \vrdg{for each of the $\tilde{k}$ isolated triangles}, and $\TP_k$  \vrdg{for each of the $t$ tadpoles} are all non-negative. So, the only term in the sum \eqref{PStarTbiasmain} that can be negative is $\C$. We will show that this term is outweighed by the other terms \vrdg{by considering three exhaustive cases}.

\medskip\noindent
\underline{Case 1}:
The graph has a tadpole, \vrdg{$t\ge 1$}. Then, using \eqref{eq:eqnrepeat}, we have
\[
(n+1)\Delta^T_{[n+1]} \geq \C+\TP_k \geq \tfrac{2k}{n} -k +k = \tfrac{2k}{n} \geq 0.
\]

\medskip\noindent
\underline{Case 2}:
The graph has no tadpoles, \vrdg{$t=0$}, but \vrdg{there exists $1\le i \le m$ such that} $k_{i} \ge 2$. Then (since there are no tadpoles) $2k \geq n$, and so $\C = \frac{2k}{n} - k \geq 1-k$. Hence, \vrdg{using \eqref{eq:eqnrepeat}, we have}
\[
(n+1)\Delta^T_{[n+1]} \geq \C+\AT_k(k_i) \geq 1-k + \tfrac{k}{3}(k_i+2) - \tfrac{2k_i}{3}\geq \tfrac{1}{3}(k-2)(k_i-1)+1 \geq 0,
\]
where the last inequality holds because $k\geq k_i \geq 2$ for a graph with a band of adjacent triangles.

\medskip\noindent
\underline{Case 3}: The graph has only isolated triangles. Then (since there are no tadpoles and no bands of adjacent triangles) $2k = n$, and so $\C = 1-k$. Hence, \vrdg{using \eqref{eq:eqnrepeat}, we have}
\[
(n+1)\Delta^T_{[n+1]} \geq \C+\IT_k = 1-k + k-1 \geq 0.
\]
\end{proof}


\subsection{Proof of Theorem \ref{thm:stargraphsglued}}

\vrdg{First note from \eqref{eq:eqnrepeat}} that, in the sum for $(n+1)\Delta^T_{[n+1]}$, the only term that can be negative is $\C$, and there is always one other term, which we henceforth call $\hat \T$, such that $\C + \hat \T\geq 0$. If the graph has a tadpole, then $\hat \T = \TP_k$. If the graph has no tadpole but a band of adjacent triangles, then $\hat \T = \AT_k(k_i)$. If the graph has neither tadpoles nor a band of adjacent triangles, then $\hat \T = \IT_k$.

\vrdg{We begin by stating a preliminary lemma that is needed to prove Theorem~\ref{thm:stargraphsglued}.}

\begin{lemma}
\label{lem:list1.5}
The partially completed star-graphs with $(n+1)\Delta^T_{[n+1]}< \tfrac32$ are
\begin{center}
\begin{tikzpicture}
\draw (1,0.3)--(2,0.3)
node[draw,isosceles triangle,isosceles triangle apex angle=60,
rotate=180, minimum size = 0.7cm, anchor=apex]{}
node[] at (2,1.2) {$G(1,1,\{\phi\})$};
\node[draw,isosceles triangle,anchor=apex,rotate = 30, isosceles triangle apex angle = 60, minimum size = 0.7cm]
(T60) at (5,0.3) {};
\node[draw, isosceles triangle,anchor=apex, rotate = 330, isosceles triangle apex angle = 60, minimum size = 0.7cm]
(T60) at (5,0.3){}
node[] at (5,1.2) {$G(0,0,\{2\})$};
\hspace{0.8cm}
\node[draw, isosceles triangle, anchor= apex, rotate = 330, isosceles triangle apex angle = 60, minimum size = 0.7cm]
(T60) at (7,0){}
node[] at (7,1.2) {$G(0,1,\{\phi\})$};
\hspace{0.8cm}
\node[draw, isosceles triangle, anchor=apex, rotate= 330, isosceles triangle apex angle = 60, minimum size=0.7cm]
(T60) at (9,0){};
\node[draw, isosceles triangle, isosceles triangle apex angle = 60, anchor=apex, rotate=210, minimum size = 0.7cm]
(T60) at (9,0){}
node[] at (9,1.2) {$G(0,2,\{\phi\})$};
\end{tikzpicture}
\end{center}
with $(n+1)\Delta^T_{[n+1]}$ equal to $\tfrac{2}{3}$, $\tfrac{2}{3}$, $0$, $1$, respectively.
\end{lemma}

\noindent
We defer the proof of Lemma~\ref{lem:list1.5} to the end of this section.

\begin{proof}[Proof of Theorem~\ref{thm:stargraphsglued}]
\vrdg{We introduce some notation to simplify the argument.} Consider two partially completed star-graphs $G_1$ and $G_2$ with $n_1$ and $n_2$ \vrdg{number of} non-center vertices, centers $C_1$ and $C_2$, $k_1$ and $k_2$ total number of triangles, respectively. Consider the new graph $G_0$ obtained by gluing $G_1$ and $G_2$ at $x \in V(G_1)\setminus C_1$ and $y\in V(G_2)\setminus C_2$. Our next lemma expresses the new total TFB of $G_0$ to be $\NB = (n_1+n_2+1)\Delta^T_{\vrdg{[n_1+n_2+1]}}$ in terms of the two old total TFB, $\OB_1 = (n_1+1)\Delta^T_{[n_1+1]}$ and $\OB_2 = (n_2+1)\Delta^T_{[n_2+1]}$, plus some extra terms. To state this lemma, we introduce \vrdg{some more} notation.
\vrdg{\begin{itemize}
\item
Let $c_1, c_2\in\{0,1,2\}$ denote the number of non-center neighbours of $x$ and $y$, respectively. Note that $c_1$ and $c_2$ are also the number of triangles that $x$ and $y$ can be part of.
\item
Let $a_1, a_2 \in\{1,2\}$ denote the number of triangles that the non-center neighbours of $x$ are part of, whenever $c_1\geq 1$ and $c_1 = 2$, respectively,
\item
Let $b_1, b_2\in\{1,2\}$ denote the same for the $c_2$ non-center neighbours of $y$ whenever these exist.
\end{itemize}}

\vrdg{We next state a key lemma.}

\begin{lemma}
\label{lem:NBcases}
Consider the gluing of two partially completed star-graphs $G_1$ and $G_2$ at $x \in V(G_1) \setminus C_1$ and $y \in V(G_2) \setminus C_2$. Let $c_1,c_2,k_1,k_2, n_1,n_2$ be defined as above, and let $k=\max\{k_1,k_2\}$ (not to be confused with the definition of $k$ in the proof of Theorem \ref{thm:stargraphs}). Then the new bias of the glued graph $\NB$ is related to the two old biases $\OB_1$ and $\OB_2$ as follows:
\begin{equation}
\label{eq:NBcases}
\NB \geq \left\{
\begin{array}{c l c c l c}
\OB_1 + \OB_2 &   & -k, &
&\mbox{if } c_1=c_2=0,\\
\OB_1+\OB_2 & + \left( \frac{c_1}{n_2} \vee \frac{c_2}{n_1}\right) & -\frac{5k}{6} & -\frac{1}{3},
&\mbox{if } c_1 \vee c_2 >0, c_1\wedge c_2 =0,\\
\OB_1+\OB_2 & +\frac{c_1}{n_2} + \frac{c_2}{n_1} &-\frac{k}{2} & -\frac{1}{3},
&\mbox {if }c_1,c_2>0.
\end{array}
\right.
\end{equation}
\end{lemma}

\noindent
We defer the proof of Lemma~\ref{lem:NBcases} to the end of the section. Note that, without loss of generality, we can let $G_1$ be the graph such that $k = \max\{k_1,k_2\} = k_1 \geq k_2$. In case $k_1 = k_2$, we let $G_1$ be the graph with more tadpoles. If the number of tadpoles are also equal, then the choice can be arbitrary.

If $k=1$, then both $G_1$ and $G_2$ have only one triangle. We leave this case as an exercise to the reader. It is easy to see that if $G_1$ has more than two tadpoles, then $\OB_1-k-\tfrac{1}{3} \geq 0$, resulting in $\NB \geq 0$. The remaining cases involving both graphs having either one or no tadpole are also easy to check. For the rest of the proof, we assume that $k\geq 2$, which implies $II + III \geq -k$ (we may use \eqref{eq:III} and \eqref{eq:II} whenever necessary). \vrdg{We split the proof into three cases.}

\medskip\noindent
\underline{Case 1}:
Consider $\OB_1$. It contains one term $\hat \T$ such that $\hat \T+\C \geq 0$. In addition, it may contain \vrdg{any of the terms in Table~\ref{tb:case1OB1}}.

\begin{table}[htbp]
\centering
\renewcommand{\arraystretch}{1.3}
\vrdg{\begin{tabular}{|l|l|}
\hline
\textbf{Term(s)} & \textbf{TFB value} \\
\hline
$\TP_k$ & $k$ \\
\hline
$\AT_k(2)$ with $k \ge 4$ & $k + \tfrac{k-4}{3}$ \\
\hline
$2\,\IT_k$ & $k + (k-2)$ \\
\hline
$\AT_k(k_i)$ for some $1 \le i \le m$ with $k_i \ge 3$ & $k + \tfrac{1}{3}((k-2)(k_i - 1) - 2)$ \\
\hline
$\AT_k(2) + \IT_k$ & $k + \tfrac{4k - 7}{3}$ \\
\hline
\end{tabular}}
\caption{\vrdg{\small Possible additional terms in $\OB_1$ other than $\hat\T$.}}
\label{tb:case1OB1}
\end{table}
Note that all the terms \vrdg{in Table~\ref{tb:case1OB1}} are greater than or equal to $k$. The resulting glued graph must have non-negative TFB because (use Lemma~\ref{lem:list1.5}) $\NB\geq \OB_1+\OB_2-k\geq \OB_1-k\geq 0$. In all other cases, either there is a term in the sum of $\OB_1$ in addition to $C +\hat \T$ that is not in the list above, or there are only two terms in the sum. We deal with these two cases below.

\medskip\noindent
\underline{Case 2}:
In addition to the term $\hat \T$ and $\C$, the sum $\OB_1$ can have additional terms not listed in Table~\ref{tb:case1OB1}. These can be:
\begin{itemize}
\item $\AT_k(2)$ with $k=2$ or $k=3$: Then we must have $\hat \T = \TP_k$, since otherwise we would have $k\geq 4$. So the graph must be
\begin{tikzpicture}
\node[draw,isosceles triangle, isosceles triangle apex angle = 60, rotate=210,anchor=apex]{};
\node[draw,isosceles triangle, isosceles triangle apex angle = 60, rotate=270,anchor=apex]{};
\draw (0,0)--(-0.45,-0.25);
\end{tikzpicture}
, denoted by $G_1 = G(1,0,\{2\})$. Hence $\OB_1 = \C+\TP_2+\AT_2(2) = \tfrac{2 \times 2}{4} - 2 + 2 + \tfrac{4}{3} =
\tfrac{7}{3}$. This, however, implies that $\NB \geq 0$, since $\NB \geq \OB_1-k \geq \tfrac{1}{3}$. Note that we cannot have $k=3$ because $AT_k(2)+IT_k$ is listed in Table~\ref{tb:case1OB1}.
\item
$\IT_k$ only: Then $\hat \T$ cannot be $TP_k$, since we have assumed $k\geq 2$. It must be either $\AT_k(k-1)$ or $\IT_2$.
\begin{itemize}
\item
$\hat \T = \AT_k(k-1)$. Here, $\OB_1 = \C +\IT_k+\AT_k(k-1) = \tfrac{2k}{k+2}-k+k-1+\tfrac{k(k+1)}{3}-\tfrac{2(k-1)}{3}$.
\begin{itemize}
\item
For $k\geq 4$, we have $\OB_1-k \geq 0$.
\item
For $k=3$, the graph $G_1$ is
\begin{tikzpicture}
\node[draw,isosceles triangle, isosceles triangle apex angle = 60, anchor=apex, rotate=90]{};
\node[draw,isosceles triangle, isosceles triangle apex angle = 60, anchor=apex, rotate=150]{};
\node[draw,isosceles triangle, isosceles triangle apex angle = 60, anchor=apex,rotate=0]{};
\end{tikzpicture},
denoted by $G(0,1,\{2\})$, and we have $\OB_1-k \geq -\tfrac{2}{15}$. If $\OB_2 \geq\tfrac{2}{15}$, then $\NB \geq 0$, since $\NB \geq \OB_1+\OB_2-k$. By Lemma \ref{lem:list1.5}, the only graph $G_2$ for which $\OB_2 \leq \tfrac{2}{15}$ is
\begin{tikzpicture}
\node[draw,isosceles triangle, isosceles triangle apex angle = 60, anchor=apex, rotate=-30]{};
\end{tikzpicture}
denoted by $G_2 = G(0,1,\{\phi\})$. In this case, since both $G_1$ and $G_2$ have no tadpoles, we have $I\geq \tfrac{c_1}{n_2}\geq\tfrac{1}{2}$. Hence $\NB \geq \OB_1+\OB_2+I-k \geq -\tfrac{2}{15}+\tfrac{1}{2} \geq 0$.
\end{itemize}
\item
$\hat \T = \IT_2$. In this case, $G_1$ is
\begin{tikzpicture}
\node[draw,isosceles triangle, isosceles triangle apex angle = 60, anchor=apex,rotate=60]{};
\node[draw,isosceles triangle, isosceles triangle apex angle = 60, anchor = apex, rotate = 150]{};
\end{tikzpicture},
denoted by $G(0,2,\{\phi\})$, and $\OB_1 = \C+\IT_2+\IT_2=1$ and $\OB_1-k= -1$.
\begin{itemize}
\item
If $\OB_2 \geq 1$, then $\NB\geq \OB_1+\OB_2-k \geq 1-1=0$.
\item
By Lemma \ref{lem:list1.5}, there are exactly three possible graphs $G_2$ with $\OB_2 < 1$:
\begin{itemize}
\item
$G_2 = G(0,1,\{\phi\})$, which is
\begin{tikzpicture}
\node[draw,isosceles triangle, isosceles triangle apex angle=60,anchor=apex, rotate=90]{};
\end{tikzpicture}
. In this case, the glued graph can only look like
\begin{tikzpicture}
\node[draw,isosceles triangle, anchor=apex,isosceles triangle apex angle = 60, rotate = 90 ](T60) at (0,0){};
\node[draw,isosceles triangle, anchor=apex,isosceles triangle apex angle=60, rotate = 90](T60) at (0.567,0){};
\node[draw,isosceles triangle, anchor=apex, isosceles triangle apex angle=60,rotate=90](T60) at (0.567*2,0){};
\end{tikzpicture}
, which has positive TFB.
\item
$G_2 = G(1,1,\{0\})$, which is
\begin{tikzpicture}
\node[draw,isosceles triangle, isosceles triangle apex angle=60,anchor=apex,rotate=210]{};
\draw (-0.45,0.35)--(0,0);
\end{tikzpicture}
, or $G(0,0,\{2\})$, which is
\begin{tikzpicture}
\node[draw,isosceles triangle, isosceles triangle apex angle = 60, rotate=60,anchor=apex]{};
\node[draw,isosceles triangle, isosceles triangle apex angle = 60, rotate = 120, anchor=apex]{};
\end{tikzpicture}.
Both these graphs have $OB_2=\tfrac{2}{3}$. Since $G_1$ has no tadpoles, $c_1\geq 1$ and so we have $I \geq
\tfrac{c_1}{n_2} \geq \frac{1}{3}$. Therefore $\NB \geq \OB_1 + \OB_2 - k + I\geq \tfrac{1}{3} +\tfrac{2}{3}-1\geq 0$.
\end{itemize}
\end{itemize}
\end{itemize}
\end{itemize}

\medskip\noindent
\underline{Case 3}:
The sum $\OB_1 = \C+\hat \T$ has no additional terms. Since $k\geq 2$, we must have $\hat \T = \AT_k(k)$. Therefore $\OB_1 = \C + \AT_k(k) = \tfrac{2k}{k+1} - k + \tfrac{k^2}{3}$.
\begin{itemize}
\item
If $k\geq 5$, then $\OB_1-k \geq 0$, resulting in $\NB\geq 0$.
\item
If $k=2,3,4$, then $\OB_1-k \geq -\tfrac32$. If $G_2$ is any graph other than those listed in Lemma \ref{lem:list1.5}, then $\NB \geq \OB_1+\OB_2-k \geq 0$, since $\OB_1-k \geq -\tfrac32$ and $\OB_2 \geq \tfrac32$. The only possibilities that remain are the ones where $\OB_2$ concerns one of the four graphs described in Lemma \ref{lem:list1.5}. Since $G_1$ has no tadpoles, one of the two occurs.
\begin{itemize}
\item
The gluing is between two non-tadpole vertices. In that case we can use the third expression of \eqref{eq:NBcases}, since for a gluing of two non-tadpole vertices we have $c_1,c_2\geq 1$. Observe that $\tfrac{c_2}{n_1} \geq \tfrac{1}{k+1}$. We have
\begin{equation*}
\begin{aligned}
\NB &\geq \OB_1 + \OB_2 - \tfrac{k}{2} - \tfrac{1}{3} + \tfrac{c_1}{n_2} + \tfrac{c_2}{n_1}\\
&\geq \tfrac{2k}{k+1} - k + \tfrac{k^2}{3} -\tfrac{k}{2} - \tfrac{1}{3} +\tfrac{1}{k+1} + \OB_2 + \tfrac{c_1}{n_2}.
\end{aligned}
\end{equation*}
For the four graphs in Lemma \ref{lem:list1.5}, we have $\OB_2+\tfrac{c_1}{n_2} \geq \tfrac{1}{2}$, resulting in $\NB\geq 0$ for all $k\in\N$.
\item
The gluing is between a tadpole and a non-tadpole vertex. This  leaves
only the case where $G_2$ is
\begin{tikzpicture}
\node[draw,isosceles triangle, isosceles triangle apex angle=60,anchor=apex,rotate=210]{};
\draw (-0.45,0.35)--(0,0);
\end{tikzpicture}
, denoted by $G(1,1,\{\phi\})$, with $\OB_2 = \frac{2}{3}$. We use the second equation of \eqref{eq:NBcases} to obtain
\begin{align*}
\NB &\geq \OB_1+\OB_2-\tfrac{5k}{6}-\tfrac{1}{3}+\tfrac{c_1}{n_2} \\
&\geq \tfrac{2k}{k+1}-k+\tfrac{k^2}{3}+\tfrac{2}{3}-\tfrac{5k}{6}-\tfrac{1}{3}+\tfrac{c_1}{3}.
\end{align*}
\begin{itemize}
\item
If either $k=4$ or $c_1=2$, then $\NB\geq 0$.
\item
Otherwise, the glued graph has $c_1=1$ and $k=2,3$, and is one of the following two graphs, both of which have positive TFB:

\begin{tikzpicture}
\node[draw,isosceles triangle, isosceles triangle apex angle=60, rotate = 0, minimum size = 0.7 cm, anchor=apex]{};
\node[draw,isosceles triangle, isosceles triangle apex angle=60, rotate = 60, minimum size = 0.7 cm, anchor=apex]{};
\node[draw,isosceles triangle, isosceles triangle apex angle=60, rotate = 120, minimum size = 0.7 cm, anchor=apex]{};
\draw (0.606217782649+0.1,-0.419)--(0.606217782649+0.8,-0.419)
node[draw,isosceles triangle, isosceles triangle apex angle=60, rotate = 180, minimum size = 0.7 cm, anchor=apex]{};
\node[draw,isosceles triangle, isosceles triangle apex angle=60, rotate = 120, minimum size = 0.7 cm, anchor=apex]
(T60) at (6,0){};
\node[draw,isosceles triangle, isosceles triangle apex angle=60, rotate = 60, minimum size = 0.7 cm, anchor=apex]
(T60) at (6,0){};
\draw (6+0.606217782649+0.1,-0.419)--(6+0.606217782649+0.8,-0.419)
node[draw,isosceles triangle, isosceles triangle apex angle = 60, rotate =180,minimum size = 0.7cm, anchor=apex]{};
\end{tikzpicture}
\end{itemize}
\end{itemize}
\end{itemize}
This exhausts all cases and completes the proof of Theorem \ref{thm:stargraphsglued}.
\end{proof}

We will now prove Lemma~\ref{lem:list1.5} and Lemma~\ref{lem:NBcases}.
\begin{proof}[Proof of Lemma~\ref{lem:list1.5}]
We know from Theorem \ref{thm:stargraphs} that the total TFB of a partially completed star-graph $G=G(t,\tilde k, \{k_i\}_{i=1}^m)$ is
\begin{equation*}
(n+1)\Delta_{[n+1]}^T = \C+\sum_{i=1}^m \AT_k(k_i) + \tilde k\, \IT_k + t\, \TP_k
\end{equation*}
with $\C$, $\AT_k$, $\IT_k$, $\TP_k$ defined in \eqref{eq:C}--\eqref{eq:TP}. As we saw in the proof of Theorem \ref{thm:stargraphs}, there is one term $\hat \T$ (either $\TP_k$, $\AT_k(k_i)$ or $\IT_k$, in this order) such that $\hat \T + \C \geq 0$. Now, in addition to $\hat \T$, if the sum contains any of the terms listed \vrdg{in Table~\ref{tb:lem1.5cases}} (which each have TFB $\geq \tfrac32$), then $(n+1)\Delta_{[n+1]}^T\geq \tfrac32$ and we are done.

\begin{table}[htbp]
\centering
\renewcommand{\arraystretch}{1.3}
\vrdg{\begin{tabular}{|l|l|}
\hline
\textbf{Term(s)} & \textbf{TFB value} \\
\hline
$\AT_k(2)$ for $k \ge 3$ & $\frac{4}{3}(k-1) \ge \frac{8}{3}$ \\
\hline
 $\IT_k$ for $k \ge 3$ & $k - 1 \ge 2$ \\
\hline
$2\,\IT_k$ & $2(k - 1) \ge 2$ \\
\hline
$\AT_k(k_i)$ for some $1 \le i \le m$, $k_i \ge 3$ & $\frac{1}{3}(k-2)(k_i+2) + \frac{4}{3} \ge 3$ \\
\hline
 $\TP_k$ for $k \ge 2$ & $k \ge 2$ \\
\hline
$2\,\TP_k$ & $k + k \ge 2$ \\
\hline
\end{tabular}}
\caption{\vrdg{\small Possible additional terms in $(n+1)\Delta^T_{[n+1]}$ other than $\hat\T$}}
\label{tb:lem1.5cases}
\end{table}
If such is not the case, then one of the following two situations occurs:

\medskip\noindent
\underline{Case 1}:
There is no other term and $(n+1)\Delta^T_{[n+1]} = \C+\hat \T$. Since $k \geq 1$, $\hat \T$ must be either $\AT_k(k)$ or $\IT_1$.
\begin{itemize}
\item
If $\hat \T = \IT_1$, then the graph is
\begin{tikzpicture}
\node[draw,isosceles triangle, isosceles triangle apex angle = 60,rotate=90]{};
\end{tikzpicture}
or $G=G(0,1,\{\phi\})$ with $(n+1)\Delta^T_{[n+1]} = 0$, as listed in the statement of the lemma.
\item
If $\hat \T= \AT_k(k)$, then the graph is $G=G(0,0,\{k\})$ with
$(n+1)\Delta^T_{[n+1]} = \C+\AT_k(k) = \frac{2k}{k+1}-k+\frac{k^2}{3}$. This
expression is at least $\tfrac{3}{2}$ for $k \geq 3$, while $k=2$
results in the graph
\begin{tikzpicture}
\node[draw,isosceles triangle, isosceles triangle apex angle = 60, rotate=60,anchor=apex]{};
\node[draw,isosceles triangle, isosceles triangle apex angle = 60, rotate = 120, anchor=apex]{};
\end{tikzpicture}
or $G(0,0,\{2\})$ with TFB $\tfrac{2}{3}$, which is listed in the statement of the lemma.
\end{itemize}

\medskip\noindent
\underline{Case 2}:
There are some terms other than $C$ and $\hat \T$, but they are not included in Table~\ref{tb:lem1.5cases}. The additional terms in the sum other than $C$ and $\hat \T$ can be
\begin{itemize}
\item
$\IT_1 + \TP_1$, in which case we must have $\hat \T = \TP_k=\TP_1$,
so the  graph is
\begin{tikzpicture}
\node[draw,isosceles triangle, isosceles triangle apex angle=60, rotate=180,anchor=apex]{};
\draw (-0.45,-0.25)--(0,0)--(-0.45,0.25);
\end{tikzpicture}
or $G(2,1,\{\phi\})$ with $(n+1)\Delta^T_{[n+1]} = \C+\TP_1+\TP_1+\IT_1 = \tfrac{2}{4}-1+1+1+0 = 1.5$.
\item
$TP_1$, in which case the graph must be the same as above.
\item
$\AT_2(2)$, in which case, since $k=2$, the term $\hat \T$ must be
$\TP_2$, so the graph is
\begin{tikzpicture}
\node[draw,isosceles triangle, isosceles triangle apex angle = 60, rotate=210,anchor=apex]{};
\node[draw,isosceles triangle, isosceles triangle apex angle = 60, rotate=270,anchor=apex]{};
\draw (0,0)--(-0.45,-0.25);
\end{tikzpicture}
or $G(1,0,\{2\})$ with $(n+1)\Delta^T_{[n+1]} = \C+\TP_2+\AT_2(2) = \tfrac{2 \times 2}{4}-2+2+\frac{4}{3} \geq \tfrac32$.
\item
$\IT_k$ with either $k=1$ or $k=2$, in which case the graph is either
$G(1,1,\{\phi\})$ given by
\begin{tikzpicture}
\node[draw,isosceles triangle, isosceles triangle apex angle=60,anchor=apex,rotate=210]{};
\draw (-0.45,0.35)--(0,0);
\end{tikzpicture} or $G(0,2,\{\phi\})$ given by
\begin{tikzpicture}
\node[draw,isosceles triangle, isosceles triangle apex angle = 60, anchor=apex,rotate=60]{};
\node[draw,isosceles triangle, isosceles triangle apex angle = 60, anchor = apex, rotate = 150]{};
\end{tikzpicture}
, both of which have a TFB smaller than $\tfrac32$, and are listed in the statement of the lemma.
\end{itemize}
This exhausts all cases and completes the proof of Lemma \ref{lem:list1.5}.
\end{proof}

We next prove Lemma~\ref{lem:NBcases}.
\begin{proof}[Proof of Lemma~\ref{lem:NBcases}]
Note that gluing $G_1$ and $G_2$ at $x$ and $y$ potentially changes the TFB of the new glued vertex (say, $xy$). It can increase the TFB of any neighbours of $x$ and $y$, including the centers of the graphs and any non-center vertices that $x$ or $y$ may share a triangle with. The gluing does not change the TFB of any vertices that are not $x$ or $y$, or any of their neighbours. We can express $\NB$ in terms of $\OB_1$ and $\OB_2$ as follows:
\begin{equation*}
\begin{aligned}
&\NB = \OB_1+\OB_2 + \left[ \frac{c_2}{n_1}+\frac{c_1}{n_2} \right] \\
&+ \left[ \mathbbm{1}_{\{c_1\geq 1\}} \frac{c_2}{a_1+1} + \mathbbm{1}_{\{c_1 = 2\}} \frac{c_2}{a_2+1}
+ \mathbbm{1}_{\{c_2\geq 1\}} \frac{c_1}{b_1+1} + \mathbbm {1}_{\{c_2 = 2\}} \frac{c_1}{b_2+1} \right]\\
&+ \left[ \frac{k_1+k_2+a_1+a_2+b_1+b_2}{c_1+c_2+2} - (c_1+c_2) - \left( \frac{k_1+a_1+a_2}{c_1+1} - c_1 \right)
- \left( \frac{k_2+b_1+b_2}{c_2+1} - c_2 \right) \right].
\end{aligned}
\end{equation*}
Here, the first term in square brackets refers to the potential increase at the centers $C_1$ and $C_2$ upon gluing, the second term in square brackets refers to the increase at any neighbours with which $x$ or $y$ share a triangle, and the third term in square brackets refers to the change in bias given by the gluing point $\Delta_{{xy}_{G_0}} - \Delta_{{x}_{G_1}}-\Delta_{{y}_{G_2}}$. In the last line, when $a_1,a_2,b_1,b_2$ are not defined, we set them equal to $0$.

On rearranging the terms above and isolating terms with $a_i$, $b_i$, $k_i$ for $i=1,2$, we get
\begin{equation}
\label{eq:NBfinalform}
\NB = \OB_1 + \OB_2 + I + II + III
\end{equation}
with
\begin{equation*}
\begin{aligned}
&I = \left[ \frac{c_2}{n_1}+\frac{c_1}{n_2} \right],\\
&II =  -\left[ k_1\frac{(c_2+1)}{(c_1+1)(c_1+c_2+2)} + k_2 \frac{(c_1+1)}{(c_2+1)(c_1+c_2+2)} \right], \\
&III =  \mathbbm{1}_{\{c_1\geq 1\}}\left[ \frac{c_2}{a_1+1} - a_1 \frac{c_2+1}{(c_1+1)(c_1+c_2+2)} \right]\\
&\qquad\qquad + \mathbbm{1}_{\{c_1 = 2\}}\left[ \frac{c_2}{a_2+1} - a_2 \frac{c_2+1}{(c_1+1)(c_1+c_2+2)} \right]\\
&\qquad\qquad + \mathbbm{1}_{\{c_2\geq 1\}}\left[ \frac{c_1}{b_1+1} - b_1 \frac{c_1+1}{(c_2+1)(c_1+c_2+2)} \right]\\
&\qquad\qquad + \mathbbm{1}_{\{c_2 = 2\}}\left[ \frac{c_1}{b_2+1} - b_2 \frac{c_1+1}{(c_2+1)(c_1+c_2+2)} \right].
\end{aligned}
\end{equation*}
Note that $I\geq 0$ and $II\leq 0$ for all glued graphs, since $k_i,c_i \geq 0$ for $i=1,2$. Furthermore, by inserting the possible values of $c_1$, $c_2$, $a_1$, $a_2$, $b_1$, $b_2$ we can easily check that
\begin{equation}
\label{eq:II}
II \geq
\begin{cases}
-k, & \mbox{ if } c_1=c_2=0,\\
-\tfrac{5}{6} k, & \mbox{ if } c_1\vee c_2 \geq 1 \mbox{ and } c_1\wedge c_2 = 0,\\
-\tfrac{1}{2} k, & \mbox{ if } c_1,c_2 \geq 1,
\end{cases}
\end{equation}
where $k=\max\{k_1,k_2\}$ and
\begin{equation}
\label{eq:III}
III \geq
\begin{cases}
0, & \mbox{ if } c_1=c_2=0, \\
-\tfrac{1}{3}, &\mbox{ if } c_1 \vee c_2 \geq 1.
\end{cases}
\end{equation}
Combining the above expressions, we get \eqref{eq:NBcases}.
\end{proof}


\section{Proofs for sparse random graphs}
\label{s.proof2}

\vrdg{In this section we prove Theorems \ref{thm:sparse}, \ref{thm:ER} and \ref{thm:dense}. The proof of Theorem~\ref{thm:sparse} follows from the definition of local weak convergence, while Theorem~\ref{thm:ER} and Theorem~\ref{thm:dense} follow from detailed combinatorial computations.}


\subsection{Proof of Theorem~\ref{thm:sparse}}

\begin{proof}
Simply note that $\Delta^T_{[n]}(G_n)  = \EE[t_{U_n}(G_n)\vrdg{\mid G_n}]$ with $U_n$ drawn uniformly at random from $[n]$. Since $t_{U_n}(G_n)$ is a local function of $G_n$ viewed from $U_n$, and $(G_n)_{n\in\N}$ \vrdg{converges locally in probability} to $(G,o)$, this implies that $\Delta^T_{[n]}(G_n)$ \vrdg{converges weakly} to $\Delta^T(G,o) = \Delta^T$ as $n\to\infty$. 
\end{proof}


\subsection{Proof of Theorem~\ref{thm:ER}}

\begin{proof}
(a) Consider ERRG\vrdg{$(n, p)$} with retention probability $p \in (0,1)$ \vrdg{on $n$ vertices}. \vrdg{Later we set $p = n^{-1} \lambda$. Note that the number of triangles at vertex $i\in [n]$ is $t_i = \frac12 \sum_{j,k\in[n]} A_{ij}A_{jk}A_{ki}$,} where the factor $\tfrac12$ is needed to avoid double counting of triangles. Since there are no self-loops, i.e., $A_{ii}=0$ for all $i \in [n]$, we need not worry about repeating indices. \vrdg{From \eqref{eq:tfpbias} and \eqref{eq:tfbgraph}, we compute the expected triangle friendship-bias as
\begin{equation}
\EE[\Delta^T_{[n]}] =\frac{1}{n} \sum_{i\in [n]} \left[ \sum_{j \in [n]} \EE\left[\frac{A_{ij}}{d_i} t_j\right]- \EE[t_i]\right].
\label{eq:ERRGproofbias}
\end{equation}}
We first compute $\EE[t_i]$ and afterwards $\EE[\frac{A_{ij}}{d_i} t_j]$. Independence gives
\begin{equation}
\EE[t_i] = \frac12\sum_{\substack{{j,k \in [n]}\\ {i\neq j\neq k}}} \mathbb{E}[A_{ij}A_{jk}A_{ki}]
= \frac12\sum_{\substack{{j,k \in [n]}\\ {i\neq j\neq k}}} p^3 = \frac12 (n-1)(n-2)p^3,
\label{eq:E2}
\end{equation}
while
\begin{equation}
\EE\left[\frac{A_{ij}}{d_i} t_j\right] = \frac12\sum_{\substack{{j,k, l \in [n]}\\ {i\neq j\neq k\neq l}}}
\mathbb{E}\left[\frac{A_{ij}}{d_{i}}A_{jk}A_{kl}A_{lj}\right]
+ \sum_{\substack{{j,k \in [n]}\\ {i\neq j\neq k}}} \mathbb{E}\left[ \frac{A_{ij}}{d_{i}}A_{jk}A_{ki} \right],
\label{eq:E1beg}
\end{equation}
\vrdg{where the first term corresponds to vertex $j$ not forming a triangle with $i$ and the second term corresponds to $j$ forming a triangle with $i$}. Note that the factor $\frac12$ does not occur in the last terms because both $j$ and $k$ are neighbours of $i$ and form a triangle with $i$. \vrdg{We compute each term in the summands of \eqref{eq:E1beg} next.} For $i,j,k,l \in [n]$ distinct,
\begin{align}
&\mathbb{E}\left[\frac{A_{ij}}{d_{i}}A_{jk}A_{kl}A_{lj}\right]
= \mathbb{E}\left[\frac{A_{ij}}{d_{i}}\right] p^3
= p^3 \mathbb{E}\left[ \left.\frac{A_{ij}}{d_{i}}\right| A_{ij} = 1\right]\mathbb{P}(A_{ij} = 1)  + 0 \nonumber\\
&= p^4 \sum_{k=0}^{n-2} \frac{1}{k+1} \begin{pmatrix} n-2 \\ k \end{pmatrix} p^k (1-p)^{n-2-k}.\label{eq:E1first}
\end{align}
For $i,j,k\in [n]$ distinct,
\begin{align}
&\mathbb{E}\left[\frac{A_{ij}}{d_{i}}A_{jk}A_{ki} \right]
= p\,\mathbb{E}\left[ \frac{A_{ij}A_{ik}}{d_{i}} \right]
= p\, \mathbb{E}\left[ \left.\frac{A_{ij}A_{ik}}{d_{i}} \right| A_{ij}A_{ik} = 1\right]\mathbb{P}(A_{ij}=1,A_{ik} =1)
+ 0 \nonumber\\
&= p^3\, \sum_{k=0}^{n-3} \frac{1}{k+2}\begin{pmatrix} n-3 \\ k \end{pmatrix} p^k (1-p)^{n-3-k}.\label{eq:E1second}
\end{align}
Substituting \eqref{eq:E1first} and \eqref{eq:E1second} into \eqref{eq:E1beg}, we have
\begin{align}
\EE\left[\frac{A_{ij}}{d_i} t_j\right]  &= \frac12 (n-1)(n-2)(n-3)p^4 \sum_{k=0}^{n-2} \frac{1}{k+1} \begin{pmatrix} n-2 \\ k \end{pmatrix} p^k (1-p)^{n-2-k}\nonumber\\
&\qquad + (n-1)(n-2)p^3 \sum_{k=0}^{n-3} \frac{1}{k+2} \begin{pmatrix} n-3 \\ k \end{pmatrix} p^k (1-p)^{n-3-k}. \label{eq:E1}
\end{align}
\vrdg{Substituting \eqref{eq:E1} and \eqref{eq:E2} into \eqref{eq:ERRGproofbias}, we obtain}
\[
\begin{aligned}
\mathbb{E}[\Delta_{[n]}^T] &= \tfrac{1}{2}(n-1)(n-2)(n-3)p^4 \sum_{k=0}^{n-2} \frac{1}{k+1}
\begin{pmatrix} n-2 \\ k \end{pmatrix} p^k (1-p)^{n-2-k}\\
&\qquad + (n-1)(n-2)p^3 \sum_{k=0}^{n-3} \frac{1}{k+2} \begin{pmatrix} n-3 \\ k \end{pmatrix} p^k (1-p)^{n-3-k}
- \tfrac{1}{2}(n-1)(n-2)p^3.
\end{aligned}
\]
\vrdg{After computing each of the above} summations, we get
\begin{align}
\mathbb{E}[\Delta_{[n]}^T] &= \tfrac{1}{2}(n-2)(n-3)p^3(1-(1-p)^{n-1}) + (n-1)p^2 - p \nonumber\\
&\qquad + p(1-p)^{n-1} - \tfrac{1}{2}(n-1)(n-2)p^3.
\label{eq:errgtfb}
\end{align}

\vrdg{Rearranging terms and setting $p = n^{-1}\lambda$, we obtain \eqref{eq:meanERRG} in Theorem~\ref{thm:ER}(a). We also note that this yields the formula in \eqref{eq:ERRG}.}

\medskip\noindent
(b) Take the limit $n \to \infty$ of \eqref{eq:meanERRG} obtained in Theorem~\ref{thm:ER}(a), \vrdg{to obtain \eqref{eq:scaledmean}}.

\medskip\noindent
(c) Since $\Delta_{[n]}^T=0$ when there are no triangles in the graph, we have
\vrdg{
\begin{equation}
\mathbb{P}(n\Delta_{[n]}^T=0) \geq \mathbb{P}(\mathrm{ERRG}(n, \tfrac{\lambda}{n}) \text{ is triangle free}).
\label{eq:trianglefree}
\end{equation}
}
The average number of triangles in $\mathrm{ERRG}(n, \frac{\lambda}{n})$ equals
\[
\binom{n}{3} \left(\frac{\lambda}{n}\right)^3.
\]
As $n\to\infty$, this average converges to $\lambda^3/3!$. Because of sparseness, the number of triangles in $\mathrm{ERRG}(n, \frac{\lambda}{n})$ has asymptotic distribution $\mathrm{POISSON}(\lambda^3/3!)$. Therefore (see \cite[Exercises 4.5 and 4.7]{vdH2017}),
\vrdg{
\[
\lim_{n\to\infty} \mathbb{P}(\mathrm{ERRG}(n, \frac{\lambda}{n}) \text{ is triangle free}) = \exp[-\lambda^3/3!].
\]
Using the above in \eqref{eq:trianglefree}, we obtain}
\[
\liminf_{n\to\infty} \mathbb{P}(n\Delta_{[n]}^T=0) \geq \exp[-\lambda^3/3!]>0.
\]
\end{proof}

\begin{remark}
\label{rem:increasing}
{\rm \vrdg{It is easy to verify that $\lambda \mapsto \zeta(\lambda)$ is strictly increasing. Indeed, write $\zeta(\lambda) = \lambda \eta(\lambda)$ with $\eta(\lambda) = \lambda-1 + (1-\tfrac12\lambda^2)\eee^{-\lambda}$. We show that $\lambda \mapsto \eta(\lambda)$ is strictly increasing on $[0, \infty)$. Observe that $\eta'(\lambda) = 1+\eee^{-\lambda}(\tfrac12\lambda^2-\lambda-1)$ and $\eta''(\lambda) = \tfrac12\eee^{-\lambda}\lambda(4-\lambda)$. Hence $\eta$ is strictly convex on $[0,4)$ and strictly concave on $(4,\infty)$. Because $\eta(0) = \eta'(0) = 0$, it follows that $\eta$ is strictly increasing on $[0,4]$. But $\eta'(\lambda) \geq 1+3\eee^{-\lambda}$ for $\lambda \geq 4$, and so $\eta$ is also strictly increasing on $[4,\infty)$.}}
\end{remark}


\subsection{Proof of Theorem~\ref{thm:CM}}

\begin{proof}
(a) \vrdg{ Consider $\cm$ with vertex set $[n]$ and degree sequence $\mathbf{d} = (d_1, \ldots, d_n)$ with $d_i >0$ for all $i\in [n]$. Let $H_i$ be the set of half edges of vertex $i\in [n]$. Define
\[
X_{ab} := \begin{cases}
1, & \text{if the distinct half edges $a$ and $b$ are paired up,}\\
0, & \text{otherwise.}
\end{cases}
\]
Then we can define the adjacency matrix of graph $\cm$ as
\[
A_{ij} = \sum_{a\in H_i}\sum_{b \in H_j} X_{ab}.
\]
Using \eqref{eq:tfpbias}, \eqref{eq:tfbgraph} and our convention (stated before Theorem~\ref{thm:CM}), we set ``the number of triangles'' at vertex $i\in [n]$ to be $t_i = \frac{1}{2}\sum_{ {j,k \in [n]} \atop {i \neq j \neq k} } A_{ij}A_{jk}A_{ki}$ and we compute the expected triangle friendship-bias as
\begin{align}
\EE[\Delta^T_{[n]}]
&=\EE\left[\frac{1}{n} \sum_{i\in [n]}\left[\frac{1}{d_i} \sum_{j \in [n]} A_{ij} t_j - t_i\right]
1_{\{d_i \neq 0\}}\right]\nonumber\\
&= \frac{1}{n} \sum_{i\in [n]} \left[\frac{1}{d_i} \sum_{{j \in [n]} \atop {j\neq i}} \EE[A_{ij} t_j]
+ \frac{1}{d_i}  \EE[A_{ii} t_i]- \EE[t_i]\right].
\label{eq:CMtfpbias}
\end{align}

We next compute $\EE[t_i]$, $\EE[A_{ii} t_i]$ and $\sum_{{j \in [n]} \atop {j\neq i}} \EE[A_{ij} t_j]$, respectively. Namely,
\begin{align}
\mathbb{E}[t_i] &= \frac{1}{2} \sum_{ {j,k \in [n]} \atop {i \neq j \neq k} }\mathbb{E}[A_{ij}A_{jk}A_{ki}] \nonumber\\
&= \frac{1}{2} \sum_{ {j,k \in [n]} \atop {i \neq j \neq k} } \sum_{{a, a' \in H_i}
\atop {a, a' \text{distinct}}}\sum_{{b, b' \in H_j} \atop {b, b' \text{distinct}}}\sum_{{c, c' \in H_k}
\atop {c, c' \text{distinct}}}\mathbb{E}[X_{ab}X_{b'c}X_{c'a'}]\nonumber\\
&= \frac{1}{2} \sum_{ {j,k \in [n]} \atop {i \neq j \neq k} }
d_i\frac{d_j}{m_1-1}\,(d_j-1)\frac{d_k}{m_1-3}\,(d_k-1)\frac{(d_i-1)}{m_1-5},
\label{eq:CMEti}
\end{align}
and, similarly,
\begin{align}
\EE[A_{ii} t_i]
&= \frac{1}{2}\sum_{{j,k \in [n]} \atop {i\neq j\neq k}} \EE[A_{ii}A_{ij}A_{jk}A_{ki}] \nonumber\\
&= \frac{1}{2}\sum_{{j,k \in [n]}\atop  {i\neq j\neq k}} d_i \frac{(d_i-1)}{m_1-1} (d_i - 2)
\frac{d_j}{m_1 -3}(d_j - 1)\frac{d_k}{m_1 - 5}(d_k -1)\frac{(d_i - 3)}{m_1 - 7}
\label{eq:CMselfloop}
\end{align}
and
\begin{align}
\sum_{ {j\in[n]} \atop {j \neq i} } \mathbb{E}[A_{ij}t_j]
&= \frac{1}{2} \sum_{ {j,k,l \in [n]} \atop {i \neq j \neq k \neq l} }\EE[A_{ij}A_{jk}A_{kl}A_{lj}]
+ \sum_{ {j,k \in [n]} \atop {i \neq j \neq k} }\EE[A^2_{ij}A_{jk}A_{ki}] \label{eq:CMEAij}\\
&= \frac{1}{2} \sum_{ {j,k,l \in [n]} \atop {i \neq j \neq k \neq l} } d_i\frac{d_j}{m_1-1}\,
(d_j-1)\frac{d_k}{m_1-3}\,(d_k-1)\frac{d_l}{m_1-5}\,(d_l-1)\frac{(d_j-2)}{m_1-7}\nonumber\\
&\qquad+ \sum_{ {j,k \in [n]} \atop {i \neq j \neq k} }
d_i\frac{d_j}{m_1-1}\,(d_j-1)\frac{d_k}{m_1-3}\,(d_k-1)\frac{(d_i-1)}{m_1-5}\nonumber\\
& \qquad+  \sum_{ {j,k \in [n]} \atop {i \neq j \neq k} } d_i \frac{d_j}{m_1 - 1}(d_i-1)\frac{(d_j-1)}{m_1-3}(d_j-2)\frac{d_k}{m_1-5}(d_k-1)\frac{(d_i-2)}{m_1-7},
\label{eq:CMAijtj}
\end{align}
}where the last term in \eqref{eq:CMEAij} has no factor $\tfrac12$ in front because the two vertices labelled $j,k$ form a triangle with $i$ and together contribute 2 to the number of triangles that the neighbours of $i$ are part of. \vrdg{Note that $A_{ij}^2 \neq A_{ij}$, since $A_{ij}$ is no longer just $0$ or $1$. Substituting \eqref{eq:CMEti}, \eqref{eq:CMselfloop}, and  \eqref{eq:CMAijtj} into \eqref{eq:CMtfpbias},} we get, with the abbreviation $\pi_k = \prod_{l=1}^k (m_1-(2l-1))$, $k \in \N$,
\[
\begin{aligned}
n\,\mathbb{E}[\Delta_{[n]}^T]
&= \sum_{i \in [n]} \Bigg[\frac{1}{2\pi_4}
\sum_{ {j,k,l \in [n]} \atop {i \neq j \neq k \neq l} } d_j(d_j-1)(d_j-2)d_k(d_k-1)d_l(d_l-1)\\
&\qquad \qquad \qquad
+ \frac{1}{\pi_3} \sum_{ {j,k \in [n]} \atop {i \neq j \neq k} } (d_i-1)d_j(d_j-1)d_k(d_k-1)\\
&\qquad \qquad \qquad\vrdg{
+ \frac{1}{\pi_4} \sum_{ {j,k \in [n]} \atop {i \neq j \neq k} } (d_i-1)(d_i -2)d_j(d_j-1)(d_j-2)d_k(d_k-1)}\\
& \qquad \qquad \qquad\vrdg{
+ \frac{1}{2\pi_4} \sum_{ {j,k \in [n]} \atop {i \neq j \neq k} } (d_i-1)(d_i -2)(d_i - 3)d_j(d_j-1)d_k(d_k-1)}\\
&\qquad \qquad \qquad
- \frac{1}{2\pi_3} \sum_{ {j,k \in [n]} \atop {i \neq j \neq k} } d_i(d_i-1)d_j(d_j-1)d_k(d_k-1)\Bigg].
\end{aligned}
\]
By working out the sums in terms of $m_k$, $1 \leq k \leq 7$, we obtain the formula in Theorem~\ref{thm:CM}(a) after a tedious but straightforward computation (available from the authors on request).

\medskip\noindent
(b) Take the limit $n \to \infty$ of the formula in Theorem \ref{thm:CM}(a) to get
\[
\begin{aligned}
\zeta(c_1,c_2,c_3) &= \frac{1}{2c_1^4}\,(c_3-3c_2+2c_1)(c_2-c_1)^2 - \frac{1}{2c_1^3}\,(c_2-3c_1+2)(c_2-c_1)^2\\
&= \frac{(c_2-c_1)^2}{2c_1^4} \big[(c_3-3c_2+2c_1) - c_1(c_2-3c_1+2)\big].
\end{aligned}
\]
The last factor in the right-hand side is non-negative because, by an easy symmetrisation argument,
\[
\frac{1}{n} \sum_{i \in [n]} d_i(d_i-1)(d_i-2)
\geq \frac{1}{n} \sum_{i \in [n]} d_i \times \frac{1}{n} \sum_{i \in [n]} (d_i-1)(d_i-2),
\]
the reason being that both $d \mapsto d$ and $d \mapsto (d-1)(d-2)$ are non-decreasing. Substitute $c_k = d^k$ to get $\zeta(c_1,c_2,c_3) = 0$ for the $d$-regular random graph.

\medskip\noindent
(c) Since $\Delta_{[n]}^T=0$ when there are no triangles in the graph, we have
\vrdg{\begin{equation}
\mathbb{P}(n\Delta_{[n]}^T=0) \geq \mathbb{P}(\cm \text{ is triangle free}).\label{eq:CMtrifree}
\end{equation}
}The latter probability can be computed in the same manner as for $\mathrm{ERRG}(n, \frac{\lambda}{n})$. Indeed, the average number of triangles in $\cm$ equals
\[
\sum_{1 \leq i < j < k \leq n} d_i\,\frac{d_j}{m_1-1}\,(d_j-1)\,\frac{d_k}{m_1-3}\,(d_k-1)\,\frac{d_i-1}{m_1-5},
\]
as is seen from the random matching of the half-edges. As $n\to\infty$, this average converges to $\nu^3/3!$ with
\[
\nu = \frac{c_2-c_1}{c_1},
\]
because the contribution of the terms where two or more indices coincide is asymptotically negligible. Because of sparseness, the number of triangles in $\cm$ has asymptotic distribution $\mathrm{POISSON}(\nu^3/3!)$.
Therefore (see \cite[Exercises 7.16 and 7.17]{vdH2017}),
\vrdg{\begin{equation*}
\lim_{n\to \infty} \mathbb{P}(\cm \text{ is triangle free}) = \exp[-\nu^3/3!].
\end{equation*}
Using the above in \eqref{eq:CMtrifree}, we obtain}
\begin{equation*}
\liminf_{n\to \infty}\mathbb{P}(n\Delta_{[n]}^T=0) \geq \exp[-\nu^3/3!] >0.
\end{equation*}
Note that $\nu=d-1$ when the graph is $d$-regular.
\end{proof}


\section{Proofs for dense random graphs}
\label{s.proof3}

\vrdg{In this section we prove Theorem \ref{thm:dense} and Theorem \ref{thm:densefail}. The proof of Theorem~\ref{thm:dense} uses the large deviation principle for graphons derived in \cite{DS2022, M2023}. The proof of Theorem~\ref{thm:densefail} is based on direct computations that use \eqref{eq:chi} and \eqref{eq:dt} from Theorem~\ref{thm:dense}.}

\begin{proof}[Proof of Theorem \ref{thm:dense}]
\vrdg{Let $\kappa \in \mathcal{W}$ be given, and let $G(n, \kappa)$ be the canonical finite graph on vertex set $[n]$ associated to $\kappa$. Recall that} the empirical graphon associated with $G_n \equiv$  \vrdg{$G(n,\kappa)$} is
\[
\kappa^{G_n}(x,y) = A_{\lceil xn \rceil,\lceil yn \rceil}
= 1_{\{\lceil xn \rceil \text{ and } \lceil yn \rceil \text{ are connected by an edge}\}},
\qquad x,y \in [0,1].
\]
\vrdg{Note that}
\[
d_i = \sum_{j \in [n]} A_{ij}, \quad t_i = \tfrac12 \sum_{j,k \in [n]} A_{ij} A_{jk} A_{ki}, \qquad i \in [n].
\]
\vrdg{It is well known that
\begin{equation}
\mbox{$\kappa^{G_n}$ converges to $\kappa$\quad $\PP$-a.s. as $n \to \infty$ in the cut-metric}
\label{eq:kappa}
\end{equation}
(see \cite[Section 2, Theorem 3.1]{AR2016}). This immediately implies that}
\begin{equation}
\lim_{n\to\infty} n^{-1} d_{\lceil xn \rceil} = \mathcal{D}(x), \quad
\lim_{n\to\infty} n^{-2} t_{\lceil xn \rceil} = \mathcal{T}(x), \quad x \in [0,1], \quad \PP\text{-a.s.},
\label{eq:dtconverge}
\end{equation}
with $\mathcal{D}(x)$ and $\mathcal{T}(x)$ as defined earlier \vrdg{in \eqref{eq:dt}}. \vrdg{Using \eqref{eq:kappa} and \eqref{eq:dtconverge} in the definition of TFB in \eqref{eq:tfpbias} and \eqref{eq:tfbgraph}, along with a careful application of the dominated convergence theorem, we get the result.}
\end{proof}


\subsection{Proof of Theorem \ref{thm:densefail}}
\label{ss.dense}

\begin{proof}
(a) For the special case where the graphon is of rank 1, i.e., $\kappa(x,y) = \nu(x)\nu(y)$ for some $\nu\colon\,[0,1] \to [0,1]$ measurable, strictly positive and continuous, we have
\[
\mathcal{D}(x) = \nu(x)\, m_1, \qquad \mathcal{T}(x) = \nu^2(x)\,m^2_2,
\]
where
\[
m_k = \int_{[0,1]} \ddd x\,\nu^k(x), \qquad k \in \N.
\]
Substitution gives
\[
\chi^T = \frac{m_3 m^2_2}{m_1} - m^3_2 = \frac{m^2_2}{m_1}\,(m_3 - m_1m_2).
\]
Rewriting, with the help of symmetrisation,
\[
\begin{aligned}
m_3 - m_1m_2 &= \int_{[0,1]} \ddd x \int_{[0,1]} \ddd y \left[\nu^3(x) - \nu(x) \nu^2(y)\right]\\
&= \tfrac12 \int_{[0,1]} \ddd x \int_{[0,1]} \ddd y \left[\nu^3(x) + \nu^3(y) - \nu(x) \nu^2(y) - \nu^2(x) \nu(y)\right].
\end{aligned}
\]
The integrand equals
\[
[\nu^2(x) - \nu^2(y)] [\nu(x)-\nu(y)] = [\nu(x) + \nu(y)] [\nu(x)-\nu(y)]^2,
\]
which is non-negative. Hence $\chi^T \geq 0$, with equality if and only if $\nu$ is constant, i.e., the graphon is that of the \emph{homogeneous Erd\H{o}s-R\'enyi random graph}.

The rank-1 property of $\kappa$ corresponds to the setting where the vertices of $G_n$ are assigned independent weights $(w_i)_{i \in [n]}$, drawn according to the probability distribution $\nu/m_1$, and vertices $i$ and $j$ are independently connected by an edge with probability $w_iw_j$. This is sometimes referred to as the \emph{Chung-Lu random graph}.

\medskip\noindent
(b) Consider the two-block graphon
\[
\kappa(x,y) =
\begin{cases}
\alpha, & \mbox{if } (x,y) \in [0,p] \times [0,p],\\
\gamma, & \mbox{if } (x,y) \in [0,p] \times [p,1] \cup [p,1] \times [0,p],\\
\beta, & \mbox{if } (x,y) \in [p,1] \times [p,1],
\end{cases}
\]
with $p \in (0,1)$ and $\alpha,\beta,\gamma \in [0,1]$, as \vrdg{summarised in} Figure \ref{fig:twoblockgraphon}.

\vspace{0.3cm}
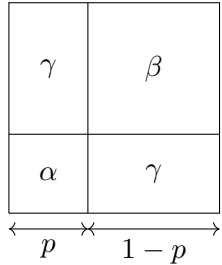
\begin{figure}[htbp]
\centering
\begin{tikzpicture}[scale=0.7]
\draw (0,0)--(0,4)--(4,4)--(4,0)--(0,0);
\draw (1.5,0)--(1.5,4);
\draw (0,1.5)--(4,1.5);
\draw [<->] (0,-0.3) -- (1.5,-0.3) node[midway,below]{$p$};
\draw [<->] (1.5,-0.3)-- (4,-0.3) node[midway,below]{$1-p$};
\node at (0.75,0.75) {$\alpha$};
\node at (2.75,2.75) {$\beta$};
\node at (0.75,2.75) {$\gamma$};
\node at (2.75,0.75) {$\gamma$};
\end{tikzpicture}
\caption{\small A two-block graphon.}
\label{fig:twoblockgraphon}
\end{figure}

\noindent
\vrdg{From \eqref{eq:dt},} the degree density is
\[
\mathcal{D}(x) = \int_{[0,1]} \kappa(x,y)\, \ddd y =
\begin{cases}
\alpha p + \gamma(1-p), &\mbox{if } x \in [0,p],\\
\gamma p + \beta(1-p), &\mbox{if } x \in [p,1],
\end{cases}
\]
and the triangle density is
\[
\begin{aligned}
\mathcal{T}(x) &= \tfrac12 \int_{[0,1]^2} \ddd y\,\ddd z\, \kappa(x,y)\kappa(y,z)\kappa(z,x)\\[0.2cm]
&=
\begin{cases}
\alpha^3 p^2 + 2\alpha\gamma^2 p(1-p) + \beta\gamma^2(1-p)^2, &\mbox{if } x \in [0,p],\\
\beta^3 (1-p)^2 + 2\beta\gamma^2 p(1-p) + \alpha\gamma^2 p^2, &\mbox{if } x \in [p,1].
\end{cases}
\end{aligned}
\]
\vrdg{Using} \eqref{eq:chi}, after a tedious but straightforward computation, $\chi^T$ is given by
\begin{equation}
\label{twoblock}
\chi^T
= \Theta_1(\alpha, \beta, \gamma, p)\cdot\Theta_2(\alpha, \beta, \gamma, p)\cdot\Theta_3(\alpha, \beta, \gamma, p)
\end{equation}
with
\[
\begin{aligned}
\Theta_1 (\alpha, \beta, \gamma, p) &= \frac{p(1-p)\gamma}{(p\gamma+(1-p)\beta)(p\alpha+(1-p)\gamma)},\\[0.2cm]
\Theta_2 (\alpha, \beta, \gamma, p)&= \big[p(\gamma-\alpha)+(1-p)(\beta-\gamma)\big],\\[0.2cm]
\Theta_3(\alpha, \beta, \gamma, p) &= \big[p^2\alpha(\gamma^2-\alpha^2)+2p(1-p)\gamma^2(\beta-\alpha)+(1-p)^2\beta(\beta^2-\gamma^2)\big].
\end{aligned}
\]
\vrdg{We will assume that $\alpha, \beta, \gamma$ are not all zero, so that $\Theta_1(\alpha, \beta, \gamma, p)$ is well defined and non-negative for all choices. Set $\alpha = 0$.
Then $\chi^T$ in \eqref{twoblock} is given by
\begin{align*}
\chi^T &= \Theta_1(0, \beta, \gamma, p)
\times \left[p\gamma- (1-p)(\gamma - \beta)\right] \times
\left[2p(1-p)\gamma^2\beta+ (1-p)^2\beta(\beta^2-\gamma^2)\right]\\
&\\
&= 2\Theta_1(0, \beta, \gamma, p) \times \left[(1-p)^2\gamma \beta
\right] \times \left[\frac{p}{1-p} -\frac{\gamma-\beta}{\gamma}\right] \times
\left[\frac{p}{1-p} - \frac{\gamma^2-\beta^2}{2\gamma^2}\right],
\end{align*}
Take any $\beta, \gamma \in (0,1)$ such that $\beta<\gamma$ and
\begin{equation}
\frac{\gamma^2-\beta^2}{2\gamma^2} < \frac{p}{1-p} < \frac{\gamma-\beta}{\gamma},
\label{eq:betagamma}
\end{equation}
to achieve $\chi^T <0$. For instance, take $p = \frac{10}{33}, \alpha = 0, \beta = \frac{1}{4} \text{ and }\gamma = \frac{1}{2}$.}
\end{proof}

\vrdg{We conclude this section with a remark containing} examples of sequences of dense random graphs for which the average TFB is either $\PP$-a.s.\ eventually positive or $\PP$-a.s.\ eventually negative.

\vrdg{\begin{remark}
\label{rem:sbm}
{\rm Here are some further examples of the $2\times 2$ stochastic block model where $\chi^T \ge 0$ and $\chi^T <0$.
\begin{itemize}
\item
\textcolor{black}{For $p=\tfrac12$, the expression in \eqref{twoblock} factorises to
\[
\chi^T = \frac{(\alpha-\beta)^2}{8(\alpha+\gamma)(\beta+\gamma)}
\left[\gamma^3+ \gamma(\alpha^2 +\alpha\beta + \beta^2)\right] \geq 0.
\]}
\item
\textcolor{black}{In \eqref{twoblock}, $\chi^T>0$  when $\beta>\gamma>\alpha$ or $\beta<\gamma<\alpha$, irrespective of the value of $p$.}
\item \textcolor{black}{If \eqref{eq:betagamma} holds and $\beta < \gamma$, it is easy to see that $\chi^T<0$ also for $\alpha$ small enough. If $p=0.4$ and $\alpha = 0.1$, $\beta = 0.3$, $\gamma = 0.8$, then $\chi^T = -8.24$, which shows that $\chi^T$ can take large negative values.}
\end{itemize}}
\end{remark}}


\bibliographystyle{unsrt}
\bibliography{references}


\end{document}